\documentclass[11pt,3p,number,sort&compress]{elsarticle}

\usepackage{lineno}
\modulolinenumbers[5]
\usepackage{moreverb}
\usepackage{subcaption}
\usepackage[colorlinks,bookmarksopen,bookmarksnumbered,citecolor=red,urlcolor=red]{hyperref}
\usepackage[american]{babel}
\usepackage{microtype}
\usepackage{float}
\usepackage{amsmath}
\usepackage{xcolor}
\usepackage[capitalize]{cleveref}
\usepackage{booktabs}
\usepackage{amssymb}
\usepackage{bm}
\usepackage{mathrsfs}
\usepackage{pgfplots}
\pgfplotsset{compat=newest}
\usepgfplotslibrary{groupplots}
\pgfkeys{/pgf/number format/.cd,1000 sep={\,}}
\usepackage[binary-units=true,exponent-product=\cdot]{siunitx}
\usepackage{url}
\usepackage{algorithm,algorithmicx,algpseudocode}
\allowdisplaybreaks[4]
\usepackage{mathdots}
\usepackage{adjustbox}

\newcommand\BibTeX{{\rmfamily B\kern-.05em \textsc{i\kern-.025em b}\kern-.08em
		T\kern-.1667em\lower.7ex\hbox{E}\kern-.125emX}}

\newcommand*{\diff}{\mathop{}\!\mathrm{d}}
\newcommand{\ue}{\mathrm{e}}
\newcommand{\im}{\mathrm{i}}

\newcommand{\Smat}{\mathbf{S}}
\newcommand{\Umat}{\mathbf{U}}
\newcommand{\Kmat}{\mathbf{K}}
\newcommand{\Vmat}{\mathbf{V}}
\newcommand{\lambdaMat}{\mathbf{\Lambda}}
\newcommand{\bigO}{\mathcal{O}}

\newcommand{\ndof}{n_{\mathrm{dof}}}
\newcommand{\Ndof}{N_{\mathrm{dof}}}
\newcommand{\Ncells}{n_{\mathrm{cells}}}
\newcommand{\nboxes}{n_{\mathrm{boxes}}}
\newcommand{\nt}{n_{\mathrm{t}}}
\newcommand{\ntNear}{n_{\mathrm{t,near}}}

\journal{Computer Methods in Applied Mechanics and Engineering}

\begin{document}
\begin{frontmatter}
\title{Fast multipole boundary element method for the acoustic analysis of finite periodic structures}

\author[address01]{Christopher Jelich}
\author[address02]{Wenchang Zhao\corref{mycorrespondingauthor}}
\cortext[mycorrespondingauthor]{Corresponding author}
\ead{Jsukya@mail.ustc.edu.cn}
\author[address02]{Haibo Chen}
\author[address01]{Steffen Marburg}

\address[address01]{Chair of Vibroacoustics of Vehicles and Machines, School of Engineering and Design, Technical University of Munich, Boltzmannstra{\ss}e 15, Garching 85748, Germany}
\address[address02]{CAS Key Laboratory of Mechanical Behavior and Design of Materials, Department of Modern Mechanics, University of Science and Technology of China, Hefei 230026, Anhui, P. R. China}

\begin{abstract}
In this work, two fast multipole boundary element formulations for the linear time-harmonic acoustic analysis of finite periodic structures are presented. Finite periodic structures consist of a bounded number of unit cell replications in one or more directions of periodicity. Such structures can be designed to efficiently control and manipulate sound waves and are referred to as acoustic metamaterials or sonic crystals. Our methods subdivide the geometry into boxes which correspond to the unit cell. A boundary element discretization is applied and interactions between well separated boxes are approximated by a fast multipole expansion. Due to the periodicity of the underlying geometry, certain operators of the expansion become block Toeplitz matrices. This allows to express matrix-vector products as circular convolutions which significantly reduces the computational effort and the overall memory requirements. The efficiency of the presented techniques is shown based on an acoustic scattering problem. In addition, a study on the design of sound barriers is presented where the performance of a wall-like sound barrier is compared to the performance of two sonic crystal sound barriers.
\end{abstract}


\end{frontmatter}

\section{Introduction}
Periodic structures are known to be very efficient in modifying the propagation of sound waves in fluids~\cite{Martinez-Sala1995, Romero-Garcia2013}. This especially holds for the sound attenuation by periodic arrangements of scatterers which are classified as sonic crystals~\cite{Martinez-Sala1995, Sigalas2005}. Sonic crystal sound barriers can lead to significant sound attenuation in certain frequency bands due to periodicity and local resonances~\cite{Maldovan2013, Melnikov2020}. Modifying the distance between the sonic crystals or changing their geometry influences the position and width of these frequency bands, which are referred to as band gaps in the case of infinite periodic arrangements. The performance of sonic crystal sound barriers is usually compared to standard wall-like sound barriers which often can be seen as periodic structures, too. Various design improvements have been proposed for wall-like sound barriers over the years, including geometric variations of the top edge~\cite{Ishizuka2004} and adding absorbing materials~\cite{Baulac2008}. Similar analyses have been carried out for sonic crystal sound barriers~\cite{Elford2011, Jean2015}, see~\cite{Fredianelli2019} for a recent review. In order to quantify the performance of designs, the acoustic behavior of the periodic structure has to be assessed within the frequency range of interest. This is described by the Helmholtz equation which can be solved by the finite element method (FEM)~\cite{Elford2011, Moheit2020}, the boundary element method (BEM)~\cite{Karimi2016, Godinho2016} and the multiple scattering theory~\cite{Amirkulova2015}.

Periodic arrangements consist of identical structures that are repeated an infinite amount of times in the directions of periodicity. Their acoustic behavior can be predicted by solving multiple acoustic eigenvalue problems of the unit cell with Floquet-Bloch boundary conditions using the FEM~\cite{Axmann1999,Morandi2016} or BEM~\cite{Gao2020}. However, this approach is not suitable for finite periodic arrangements, i.e.,\ arrangements where an identical structure is repeated only a finite amount of times. It neglects possible scattering effects at the edges of the finite periodic arrangement as well as reflections from the ground. Other modeling approaches truncate the geometry such that only a sufficiently large section of the periodic arrangement is analyzed. \citet{Reiter2017} proposed a periodic FEM which applies periodic boundary conditions to the unit cell. Analogously, \citet{Lam1999} introduced a periodic boundary element formulation which includes an infinite sum of Green's functions. Truncating this sum yields the quasi-periodic BEM which is capable of sufficiently approximating the acoustic behavior within a unit cell~\cite{Fard2015, Fard2017}. A convergence study on the truncation number is reported by~\citet{Jean2015}. However, applying periodic boundary conditions is not suitable when only the geometry is periodic but the loading and the solution are assumed to be aperiodic. A third modeling approach is to assume an infinite extent of the periodic arrangement in one spatial direction. This leads to two-dimensional numerical models which neglect edge effects along the third dimension. \Citet{Cavalieri2019} employed the FEM to analyze an infinitely long sound barrier with a constant cross section in which Helmholtz resonators and quarter-wavelength resonators are horizontally placed. With respect to the BEM, \citet{Duhamel1996} extended the two-dimensional modeling approach by incoherent line sources and refers to it as 2.5-dimensional boundary element method. In contrast, \citet{Chalmers2009} and \citet{Elford2011} assumed that the height of the sonic crystals is infinite and that ground effects are negligible. This allows to modify the cross section of the sound barrier along its length. Whenever the aforementioned assumptions are not valid, a three-dimensional acoustic analysis has to be carried out. In this regard, two major characteristics of the BEM favor its application: The BEM implicitly fulfills the Sommerfeld radiation condition and reduces the problem's dimension by one, such that only the sound radiating surface has to be discretized instead of the surrounding acoustic volume.

A recent attempt to reduce the computational effort of boundary element analyses of finite periodic arrangements was presented by~\citet{Karimi2016}. They identified that the boundary element discretization of a Helmholtz problem with periodic geometry yields a block Toeplitz system of linear equations. Utilizing this special matrix structure significantly accelerates the setup and solution times which enables three-dimensional analyses of various types of small periodic structures on todays' standard desktop computers~\cite{Karimi2017,Jelich2021}. Despite the significant reduction on the memory requirements, analyzing moderately large periodic structures with many degrees of freedom (dofs) is still infeasible. For periodicity in one direction, the storage of $2M_x -1$ dense matrices of size $\ndof\times{}\ndof$ is required, with $M_x$ denoting the number of unit cells in the $x$-direction and $\ndof$ denoting the number of degrees of freedom of a unit cell. Considering periodicity in additional directions increases the storage costs even further which scale of order~$\bigO(\Ncells\ndof^2)$ with the total number of unit cells~$\Ncells$. Well-known techniques to reduce the computational complexity of boundary element analyses in both time and memory are fast boundary element formulations. They allow to express the dense system matrix as a function of sparse matrices. The most common approaches are the fast multipole method (FMM)~\cite{Greengard1987, Nishimura2002} and hierarchical matrices~\cite{Boerm2003,Hackbusch2002}. Both are applicable to problems with arbitrary geometry but can be optimized for finite and infinite periodic arrangements. Analyzing the latter with the boundary element method involves an infinite sum of the Green's function also denoted as periodic Green's function, periodic sum or infinite lattice sum. This sum is usually truncated after a certain number of terms and evaluated by employing the FMM~\cite{Rokhlin1994,Challacombe1997,Gumerov2014}. \citet{Otani2007} introduced an FMM for the two-dimensional Helmholtz equation with periodic geometries which is extended to the three-dimensional case in~\cite{Niino2012,Ziegelwanger2017}. Similar approaches based on multipole expansions of the periodic Green's function exist for the method of fundamental solutions (MFS) which is closely related to the BEM~\cite{Liu2016}. Furthermore, \citet{Yan2018} introduced a kernel independent FMM for periodic Laplace and Stokes problems. In all approaches, the geometry is of infinite extent and the solution as well as boundary conditions and incident wave fields are assumed to be quasi-periodic. In the case of sound barriers, this assumption is often violated due to the limited height and width or aperiodic incident wave fields. Hence, finite periodic arrangements need to be analyzed. \citet{Amado-Mendes2019} applied the hierarchical matrix BEM to a finite periodic array of acoustic scatterers. They reported savings in the setup time and memory requirements due to the underlying Toeplitz structure but did not accelerate the matrix-vector products. \citet{Gumerov2005} combined a T-matrix based approach and the FMM to accelerate the solution of acoustic problems that feature arrangements of arbitrarily shaped scatterers.

We propose two novel fast multipole boundary element formulations for acoustic problems with finite periodic geometries. They exhibit low memory requirements, fast matrix-vector multiplication and allow to efficiently analyze periodic structures of finite extent. The methods do not require any periodicity of the boundary conditions, incident fields or solutions. In contrast to the aforementioned approaches, our techniques rigorously make use of the multilevel block Toeplitz matrix structure which occurs when discretizing periodic structures with the BEM and applying the FMM. Identifying that the fast multipole operators are block Toeplitz matrices is a key point of the algorithms and drives its performance. The remainder of the paper is outlined as follows. \Cref{sec:BEMformulations} presents the periodic boundary element formulation of~\cite{Karimi2016} and extends the formulation to specific symmetry problems. \Cref{sec:toeplitzFMBEM} outlines the low frequency fast multipole method and introduces our fast multipole periodic boundary element formulations. The formulations are validated numerically in terms of an acoustic scattering problem and a sound barrier design study in~\cref{sec:numericalExamples}.

\section{Boundary element method for finite periodic arrays} \label{sec:BEMformulations}
\subsection{Boundary element formulation} \label{sec:conventionalBEM}
Consider the acoustic Helmholtz equation
\begin{align}
	\label{eq:helmholtz}
	\nabla^2 p(\mathbf{x})+k^2 p(\mathbf{x}) = 0, \quad \forall \mathbf{x} \in \Omega 
  \text{ ,}
\end{align}
where $\nabla^2$ is the Laplace operator, $p$ is the acoustic pressure, $k=\omega/c$ is the wave number, $\omega = 2\pi f$ is the angular frequency with frequency~$f$ and $c$ is the speed of sound of the acoustic medium within the domain $\Omega$. The harmonic time dependence $\ue^{-\im \omega t}$ with the imaginary unit $\im$ is omitted throughout this paper. The fluid particle velocity~$v_\mathrm{f}$ on the boundary $\Gamma$ of the acoustic domain is related to the sound pressure by
\begin{align}
	\dfrac{\partial p(\mathbf{x})}{\partial \mathbf{n}(\mathbf{x})} = \im \omega \rho v_\mathrm{f}(\mathbf{x}) \text{ ,}
\end{align}
where $\mathbf{n}(\mathbf{x})$ is the outward pointing normal vector at~$\mathbf{x}\in\Gamma$ and $\rho$ is the fluid density. A well-posed problem is obtained by introducing the admittance boundary condition
\begin{align}
	\label{eq:admittance}
  \im \omega \rho v_\mathrm{f}(\mathbf{x}) = \im k \beta(\mathbf{x}) p(\mathbf{x}) \text{ ,}
\end{align}
where $\beta$ is the normalized surface admittance which can be expressed by the surface admittance~$Y$ through~\mbox{$Y(\mathbf{x}) = \beta(\mathbf{x}) / \rho c$}.

Applying Green's second theorem yields the conventional boundary integral equation referred to as CBIE~\cite{Marburg2018}
\begin{align}
	\label{eq:cbie}
	c(\mathbf{x}) p(\mathbf{x}) 
	+ \int_{\Gamma} \dfrac{\partial G(\mathbf{x}, \mathbf{y})}{\partial \mathbf{n}(\mathbf{y})} p(\mathbf{y}) \diff \Gamma(\mathbf{y}) 
	= \int_{\Gamma}G(\mathbf{x}, \mathbf{y}) \dfrac{\partial p(\mathbf{y})}{\partial \mathbf{n}(\mathbf{y})} \diff \Gamma(\mathbf{y})
	+ p^\mathrm{inc}(\mathbf{x}) \text{ ,}
\end{align}
where $\mathbf{x}$ denotes the field point and $\mathbf{y}$ denotes the source point. The solid angle $c(\mathbf{x})$ equals $1/2$ if the boundary around~$\mathbf{x}$ is smooth and equals $1$ if $\mathbf{x}\in\Omega$. Further, $G(\mathbf{x}, \mathbf{y})$ denotes the Green's function and $p^{\mathrm{inc}}(\mathbf{x})$ is the incident acoustic pressure field. The Green's functions of 3D full-space and half-space acoustic problems are given as
\begin{align}
	\label{eq:green3dFull}
	G(\mathbf{x},\mathbf{y}) &= 
		\dfrac{1}{4\pi} \dfrac{\ue^{\im k |\mathbf{x}-\mathbf{y}|}} {|\mathbf{x}-\mathbf{y}|} 
		\text{ ,} & \quad &\mathbf{x}, \mathbf{y} \in \mathbb{R}^3 
    \text{ ,} \\
  \label{eq:green3dHalf}
  G_\mathrm{h}(\mathbf{x},\mathbf{y}) &= 
    \dfrac{1}{4\pi} \dfrac{\ue^{\im k |\mathbf{x}-\mathbf{y}|}} {|\mathbf{x}-\mathbf{y}|} 
    + R_{\mathrm{p}} \dfrac{1}{4\pi} \dfrac{\ue^{\im k |\mathbf{x}-\hat{\mathbf{y}}|}} {|\mathbf{x}-\hat{\mathbf{y}}|} 
    \text{ ,} & \quad &\mathbf{x}, \mathbf{y}, \hat{\mathbf{y}} \in \mathbb{R}^3
    \text{ ,}
\end{align}
respectively, where $|\cdot|$ is the $l^2$-norm. Further, $R_{\mathrm{p}}$ is the reflection coefficient and $\hat{\mathbf{y}}$ is the mirror image of the source point~$\mathbf{y}$ with respect to the plane that divides both half-spaces. Taking the derivative of \cref{eq:cbie} with respect to the outward normal vector and assuming that the boundary around $\mathbf{x}$ is smooth yields the hypersingular boundary integral equation (HBIE)
\begin{align}
	\label{eq:hbie}
	\dfrac{1}{2} \dfrac{\partial p(\mathbf{x})}{\partial \mathbf{n}(\mathbf{x})} 
	+ \int_{\Gamma} \dfrac{\partial^2 G(\mathbf{x}, \mathbf{y})}{\partial \mathbf{n}(\mathbf{x}) \partial \mathbf{n}(\mathbf{y})} p(\mathbf{y}) \diff \Gamma(\mathbf{y}) 
	= \int_{\Gamma}\dfrac{\partial G(\mathbf{x}, \mathbf{y})}{\partial \mathbf{n}(\mathbf{x})} \dfrac{\partial p(\mathbf{y})}{\partial \mathbf{n}(\mathbf{y})} \diff \Gamma(\mathbf{y}) 
	+ \dfrac{\partial p^\mathrm{inc}(\mathbf{x})}{\partial \mathbf{n}(\mathbf{x})} \text{ .}
\end{align}
The linear combination of the CBIE and HBIE results in the Burton and Miller formulation, i.e.,\ 
\begin{align}
	\label{eq:bm}
	\begin{aligned}
	\dfrac{1}{2} p(\mathbf{x}) &+ \int_{\Gamma}\left[\dfrac{\partial G(\mathbf{x}, \mathbf{y})}{\partial \mathbf{n}(\mathbf{y})} 
	+ \alpha \dfrac{\partial^2 G(\mathbf{x}, \mathbf{y})}{\partial \mathbf{n}(\mathbf{x}) \partial \mathbf{n}(\mathbf{y})}\right] p(\mathbf{y}) \diff \Gamma(\mathbf{y}) 
  + \dfrac{\alpha}{2}\dfrac{\partial p(\mathbf{x})}{\partial \mathbf{n}(\mathbf{x})} \\
	&- \int_{\Gamma}\left[G(\mathbf{x}, \mathbf{y}) 
	+ \alpha \dfrac{\partial G(\mathbf{x}, \mathbf{y})}{\partial \mathbf{n}(\mathbf{x})}\right] \dfrac{\partial p(\mathbf{y})}{\partial \mathbf{n}(\mathbf{y})} \diff \Gamma(\mathbf{y}) 
	= p^\mathrm{inc}(\mathbf{x}) + \alpha \dfrac{\partial p^\mathrm{inc}(\mathbf{x})}{\partial \mathbf{ n}(\mathbf{x})}
	\text{ ,}
	\end{aligned}
\end{align}
which yields unique solutions at all frequencies~\cite{Burton1971}. The coupling parameter~$\alpha$ is a complex-valued scalar with a non-vanishing imaginary part, i.e., $\text{Im}(\alpha)\neq 0$. A value of~$-\im/k$ is optimal in the present case~\cite{Marburg2016}.

The collocation boundary element method is applied to discretize~\cref{eq:bm}. Quadrilateral boundary elements with second order, $C^0$-continuous Lagrange polynomials approximate the geometry, whereas discontinuous Lagrange polynomials of variable order approximate the physical quantities, i.e.,\ acoustic pressure, fluid particle velocity and boundary admittance. This leads to the following system of linear equations
\begin{align}
	\label{eq:matrixBEM}
	\left(\mathbf{H}-\mathbf{GY}\right) \mathbf{p} = \mathbf{p}^\mathrm{inc} \text{ ,}
\end{align}
where $\mathbf{H}$ and $\mathbf{G}$ are the dense boundary element coefficient matrices, $\mathbf{Y}$ is the block diagonal admittance matrix and the vectors $\mathbf{p}^\mathrm{inc}$ and $\mathbf{p}$ store the incident sound pressure values and unknown sound pressure values at the collocation points, respectively. The collocation boundary element method is described in great detail in~\cite{Marburg2018,Wu2000}.

\subsection{BEM for finite periodic structures} \label{sec:toeplitzBEM}
Both the assembly and storage of the boundary element matrices become infeasible for medium to large-scale acoustic problems. A remedy is found by~\citet{Karimi2016} for problems which feature a finite periodic geometry. In such cases, the system matrix~$(\mathbf{H} - \mathbf{G}\mathbf{Y})$ is a multilevel block Toeplitz matrix due to the translation invariance of the Green's function, c.f.~\cref{eq:green3dFull}. Block Toeplitz matrices are a special class of matrices which have constant blocks along each diagonal. Within this context, the term multilevel implies that the matrix blocks are block Toeplitz matrices itself. Representing the system matrix as a multilevel block Toeplitz matrix significantly reduces the storage cost.

\begin{figure}
	\centering
	\includegraphics[width=0.35\textwidth]{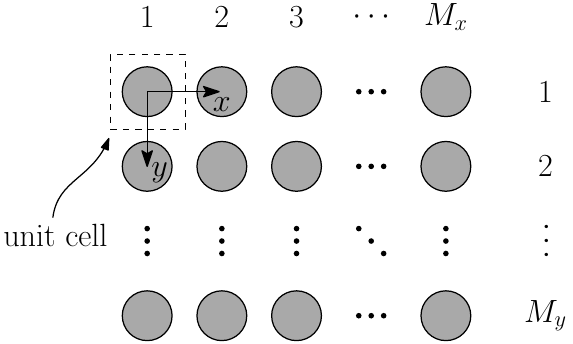}
	\caption{A configuration of a finite periodic array of scatters with periodicity in the $x$- and $y$-direction. The number of periodic segments in each direction is~$M_x$ and $M_y$. The geometry of each scatterer coincides with the geometry of the unit cell.}
	\label{fig:PeriodicArray2d}
\end{figure}%
Consider the finite periodic array of scatterers in~\cref{fig:PeriodicArray2d}. The scatterers are arranged in a regular pattern with~$M_x$ unit cells in the $x$-direction and $M_y$ unit cells in the $y$-direction. Applying a boundary element discretization yields
\begin{align}
  \label{eq:toeplitzBEM}
  \mathbf{T} \mathbf{p} = \mathbf{p}^\mathrm{inc} \text{ .}
\end{align}
where~$\mathbf{T}$ is a $2$-level block Toeplitz matrix and an exact representation of~$(\mathbf{H} - \mathbf{G}\mathbf{Y})$. It has a size of \mbox{$M_x M_y \ndof \times M_x M_y \ndof$}, with $\ndof$ denoting the number of degrees of freedom of a unit cell and reads
\begin{align}
  \label{eq:blockToeplizMatrix}
  \mathbf{T} = 
  \begin{bmatrix}
    \mathbf{T}_{0}^1   & \mathbf{T}_{-1}^1  & \cdots  				  & \cdots 
      & \mathbf{T}_{1-M_y}^1 \\
    \mathbf{T}_{1}^1   & \mathbf{T}_{0}^1   & \mathbf{T}_{-1}^1 & \cdots 
      & \mathbf{T}_{2-M_y}^1 \\
    \vdots 		     	   & \mathbf{T}_{1}^1   & \ddots 					& \ddots 
      & \vdots \\
    \vdots 		     	   & \vdots				    & \ddots 					& \ddots 
      & \mathbf{T}_{-1}^1 \\
    \mathbf{T}_{M_y-1}^1 & \mathbf{T}_{M_y-2}^1	&	\cdots	& \mathbf{T}_{1}^1
      & \mathbf{T}_{0}^1  \\
  \end{bmatrix}
  \text{ .}
\end{align}
The individual entries are block Toeplitz matrices $\mathbf{T}_{j}^1$, $j = 1-M_y {,} \, 2-M_y {,} \, \ldots {,} \, 0 {,} \, \ldots {,} \, M_y-1$ itself with $(\cdot)^1$ denoting the first level of periodicity. The matrix blocks have a size of \mbox{$M_x \ndof \times M_x \ndof$} and read
\begin{align}
  \label{eq:blockToeplizMatrixLevel2}
  \mathbf{T}_{j}^1 =  
  \begin{bmatrix}
    \mathbf{T}_{0}^2   & \mathbf{T}_{-1}^2  & \cdots  				  & \cdots 
      & \mathbf{T}_{1-M_x}^2 \\
    \mathbf{T}_{1}^2   & \mathbf{T}_{0}^2   & \mathbf{T}_{-1}^2 & \cdots 
      & \mathbf{T}_{2-M_x}^2 \\
    \vdots 		     	   & \mathbf{T}_{1}^2   & \ddots 					& \ddots 
      & \vdots \\
    \vdots 		     	   & \vdots				    & \ddots 					& \ddots 
      & \mathbf{T}_{-1}^2 \\
    \mathbf{T}_{M_x-1}^2 & \mathbf{T}_{M_x-2}^2	&	\cdots	& \mathbf{T}_{1}^2
      & \mathbf{T}_{0}^2  \\
  \end{bmatrix}_{j}
  \text{ .}
\end{align}
They store dense $\ndof\times \ndof$ matrices~$\mathbf{T}_{i}^2$, $i = 1-M_x {,} \, 2-M_x {,} \, \ldots {,} \, 0 {,} \, \ldots {,} \, M_x-1$, which are defined by the corresponding entries of the discretized boundary integral equation~\eqref{eq:bm}.

Every multilevel block Toeplitz matrix is uniquely defined by its first block row and column in each level. Therefore, the periodic boundary element formulation requires the storage of~$(2M_x - 1)(2M_y - 1)\ndof^2$ matrix entries. In contrast, the conventional boundary element method demands the storage of~$(M_x M_y \ndof)^2 = \Ndof^2$ matrix entries where~$\Ndof$ denotes the total number of dofs. The formulation will be referred to as periodic boundary element method (PBEM) throughout this work. It is referred to~\cite{Karimi2016} for a more detailed derivation and to~\cite{Karimi2017} for an extension to an arbitrary number of periodic directions and certain kinds of half-space problems. In either case, the storage cost scales of order~$\bigO(\Ncells\ndof^2)$ with~$\Ncells$ denoting the number of periodic unit cells.

\subsection{BEM for finite periodic structures in half-spaces} \label{sec:CBEM_halfspace}
\begin{figure}
	\centering
	\includegraphics[width=0.38\textwidth]{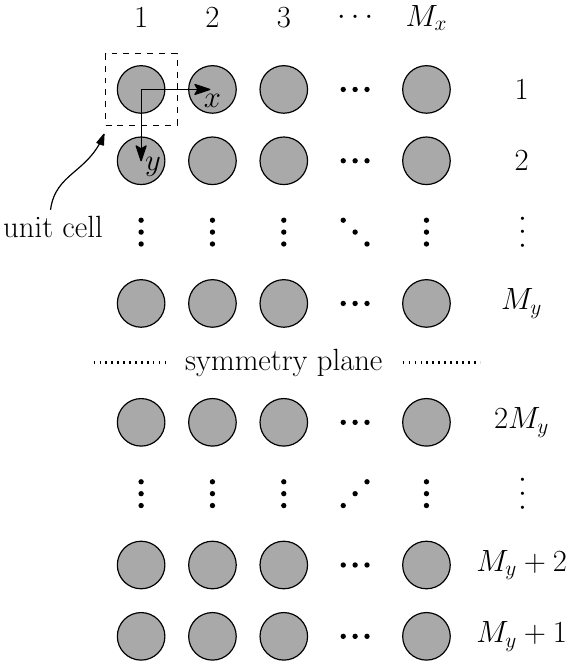}
	\caption{A configuration of a finite periodic array of scatters with periodicity in the $x$- and $y$-direction as well as symmetry along the $xz$-plane. The number of periodic segments in each direction is~$M_x$ and $2M_y$. The geometry of each scatterer coincides with the geometry of the unit cell.}
	\label{fig:PeriodicArray2dSymmetry}
\end{figure}%
The system matrix $\mathbf{H} - \mathbf{GY}$ exhibits a block Toeplitz structure for periodic geometries in half-spaces as long as the half-space Green's function \eqref{eq:green3dHalf} features translation invariance. This is the case whenever the directions of periodicity are parallel to the symmetry plane~\cite{Karimi2017}. A periodic geometry for which this does not hold is shown in~\cref{fig:PeriodicArray2dSymmetry}. The scatterers' location is periodic in the $x$ and $y$-direction and symmetric with respect to the $xz$-plane. A straightforward boundary element discretization of this problem does not yield a block Toeplitz system matrix since the periodicity in the~$y$-axis is perpendicular to the~$xz$-plane. This work introduces a remedy by taking the additive nature of the half-space Green's function~\eqref{eq:green3dHalf} into account. The idea is to split the integrals of~\cref{eq:bm} into integrals including the first summand of the half-space Green's function and integrals including the second summand. This yields two separate matrices subsequent to applying the boundary conditions and the collocation boundary element discretization. The corresponding linear system of equations reads
\begin{equation}
  \label{eq:toeplitzBEMsymmetry}
  \big(\underbrace{(\mathbf{H} - \mathbf{GY})}_{\mathbf{T}} + 
        \underbrace{(\hat{\mathbf{H}} - \hat{\mathbf{G}}\mathbf{Y})}_{\hat{\mathbf{T}}} \big)
    \mathbf{p} = \mathbf{p}^{\mathrm{inc}}
  \text{ .}
\end{equation}
The matrices~$\mathbf{H}$ and $\mathbf{G}$ stem from the first summand of the half-space Green's function which corresponds to the full-space Green's function. Hence, the multilevel block Toeplitz matrix~$\mathbf{T}$ of~\cref{eq:toeplitzBEMsymmetry} coincides with the system matrix of the corresponding full-space problem in~\cref{eq:toeplitzBEM}. The matrices $\hat{\mathbf{H}}$ and $\hat{\mathbf{G}}$ stem from the second summand of the half-space Green's function and build the matrix~$\hat{\mathbf{T}}$ which, in contrast to~$\mathbf{T}$, has constant blocks along its anti-diagonals, i.e.,\ 
\begin{align}
  \label{eq:blockToeplizMatrixSymmetry}
  \hat{\mathbf{T}} =
    \begin{bmatrix}
    \hat{\mathbf{T}}_{0}^1   & \hat{\mathbf{T}}_{-1}^1  & \cdots  	
    & \hat{\mathbf{T}}_{-M_y+1}^1 	& \hat{\mathbf{T}}_{-M_y}^1 \\
    \hat{\mathbf{T}}_{-1}^1   & \hat{\mathbf{T}}_{-2}^1   & \iddots 
    & \iddots 	& \hat{\mathbf{T}}_{-M_y-1}^1 \\
    \vdots 	& \iddots   & \iddots 
    & \iddots  & \vdots \\
    \hat{\mathbf{T}}_{-M_y+1}^1	& \iddots 	& \iddots 					
    & \hat{\mathbf{T}}_{-2My+2}^1 	& \hat{\mathbf{T}}_{-2M_y+1}^1 \\
    \hat{\mathbf{T}}_{-M_y}^1 & \hat{\mathbf{T}}_{-M_y-1}^1	&	\cdots
    & \hat{\mathbf{T}}_{-2M_y+1}^1 	& \hat{\mathbf{T}}_{-2My}^1  \\
    \end{bmatrix}
  \text{ .}
\end{align}
Structured matrices in the form of~\cref{eq:blockToeplizMatrixSymmetry} are referred to as multilevel block Hankel matrices. Since adding~$\mathbf{T}$ and~$\hat{\mathbf{T}}$ leads to an unstructured matrix, the benefits of efficient memory storage as well as efficient matrix-vector products are only preserved when assembling and storing both matrices separately. The extension of the PBEM to this special type of symmetry condition therefore requires twice the storage, i.e.,\ \mbox{$2(2M_x - 1)(2M_y - 1)\ndof^2$} matrix entries in the present case.

\subsection{Solution scheme for the periodic boundary element method} \label{sec:PBEM_solution}
Both direct and iterative solvers are available for the solution of the multilevel block Toeplitz system~\eqref{eq:toeplitzBEM}. However, an efficient direct solution of~\cref{eq:toeplitzBEMsymmetry} is infeasible since the system matrix is a sum of two structured matrices. Therefore, the focus is set on iterative solvers~\cite{Jin2003,Chan2007}.

The computational effort of iterative solvers is driven by the efficiency of the matrix-vector product. Multiplying a multilevel block Toeplitz matrix by a vector scales quasi-linearly in time when expressed by circular convolutions~\cite{Golub1989}. The present block Toeplitz matrix~$\mathbf{T}$ features dense matrices on its lowest level which leads to an asymptotic complexity of~$\bigO(\ndof^2\Ncells\mathrm{log}(\Ncells))$. In order to achieve this complexity, the multilevel block Toeplitz matrix is embedded into a multilevel block circulant matrix. This is a special type of multilevel block Toeplitz matrix where each block row is a rightward circular shift of the first block row~\cite{Davis1994,Gray2005}. For a $2$-level block Toeplitz matrix~$\mathbf{T}$ of the problem in~\cref{fig:PeriodicArray2d} with periodicity in two directions, the corresponding $2$-level block circulant matrix is defined as
\begin{equation}
  \label{eq:BlockCirculantMatrix}
  \mathbf{C} =
  \begin{bmatrix}
    \mathbf{C}^{1}_{0}  & \mathbf{C}^{1}_{-1} & \cdots 			 & \mathbf{C}^{1}_{2} &
    \mathbf{C}^{1}_{1}  \\
    \mathbf{C}^{1}_{1}  & \mathbf{C}^{1}_{0}	& \mathbf{C}^{1}_{-1} & \cdots      & 
    \mathbf{C}^{1}_{2}  \\
    \vdots 		   & \mathbf{C}^{1}_{1}	& \ddots  		 & \ddots      & \vdots       \\
    \mathbf{C}^{1}_{-2} & \vdots				& \ddots       & \ddots 		 & 
    \mathbf{C}^{1}_{-1} \\
    \mathbf{C}^{1}_{-1} & \mathbf{C}^{1}_{-2} & \cdots       & \mathbf{C}^{1}_{1} &
    \mathbf{C}^{1}_{0}
  \end{bmatrix}
  \text{ .}
\end{equation}
It stores $1$-level block circulant matrices~$\mathbf{C}_{j}^1$, which reads
\begin{align}
  \mathbf{C}_{j}^1 =  
  \begin{bmatrix}
    \mathbf{T}_{0}^2   & \mathbf{T}_{-1}^2  & \cdots  				  & \mathbf{T}_{2}^2	&
    \mathbf{T}_{1}^2 \\
    \mathbf{T}_{1}^2   & \mathbf{T}_{0}^2   & \mathbf{T}_{-1}^2 & \cdots 	&
    \mathbf{T}_{2}^2 \\
    \vdots 		     	   & \mathbf{T}_{1}^2   & \ddots 					& \ddots 		& 
    \vdots \\
    \mathbf{T}_{-2}^2	 & \vdots				    & \ddots 					& \ddots 			&
    \mathbf{T}_{-1}^2 \\
    \mathbf{T}_{-1}^2 & \mathbf{T}_{-2}^2	&	\cdots	& \mathbf{T}_{1}^2 		&
    \mathbf{T}_{0}^2  \\
  \end{bmatrix}_{j}
  \text{ .}
\end{align}
The first block row of~$\mathbf{C}_j^1$ consists of the unique matrix blocks of~$\mathbf{T}_j^1$ defined in~\cref{eq:blockToeplizMatrixLevel2}. The blocks are concatenated into the matrix
\begin{equation}
  \mathbf{Q} =
  \begin{pmatrix}
    \mathbf{T}^2_{0}  & \mathbf{T}^2_{-1} & \cdots & \mathbf{T}^2_{1-M_x} & \mathbf{T}^2_{M_x-1} & \cdots & \mathbf{T}^2_{2} & \mathbf{T}^2_{1} 
  \end{pmatrix}
  \text{ ,}
\end{equation}
and shifted in a rightward direction to build the block rows of~$\mathbf{C}_j^1$. Concatenating the unique entries of~$\mathbf{C}_j^1$ and introducing rightward shifts for each block row then builds the multilevel block circulant matrix~$\mathbf{C}$ as in~\cref{eq:BlockCirculantMatrix}. This matrix can be block diagonalized by applying Fourier transformations, i.e.,\ 
\begin{align}
  \label{eq:diagonalized}
  \mathbf{C} = \mathfrak{F}^{-1} \lambdaMat \mathfrak{F} 
  \text{ .}
\end{align}
The Fourier transform and inverse Fourier transform are defined as
\begin{align}
  &\mathfrak{F} = \mathbf{F}_{2M_y-1} \otimes \mathbf{F}_{2M_x-1} \otimes 
  \mathbf{I}_{\ndof}
  \text{ ,} \\
  &\mathfrak{F}^{-1} = \mathbf{F}_{2M_y-1}^{-1} \otimes \mathbf{F}_{2M_x-1}^{-1} \otimes
  \mathbf{I}_{\ndof}
  \text{ ,}
\end{align}
where~$\otimes$ denotes the Kronecker product and $\mathbf{I}_{\ndof}$ denotes the $\ndof \times \ndof$ identity matrix. The Fourier matrix as well as its inverse are given by
\begin{align}
  &\mathbf{F}_{m} = 
  \left(\left(\mathrm{e}^{-\im 2\pi / m}\right)^{ij}\right)_{i,j = 0,\,\ldots,\,m-1}
  \text{ ,} \\
  &\mathbf{F}_{m}^{-1} = \dfrac{1}{m} \mathbf{F}_{m}^{*}
  \text{ ,}    
\end{align}
with the complex conjugate transpose~$(\cdot)^{*}$~\cite{Golub1989}. The block diagonal matrix~$\lambdaMat$ of~\cref{eq:diagonalized} reads
\begin{equation} 
  \label{eq:blockDiagonal}
  \mathbf{\lambdaMat} = \mathrm{diag} 
  \left(
    \mathfrak{C}_{0}^{1},\ \mathfrak{C}_{1}^{1},\ \dots,\ \mathfrak{C}_{-2}^{1},\ \mathfrak{C}_{-1}^{1}
  \right)
  \text{ ,}
\end{equation}
with $[\mathfrak{C}_{0}^{1}\ \mathfrak{C}_{1}^{1}\ \dots\ \mathfrak{C}_{-2}^{1}\ \mathfrak{C}_{-1}^{1}] = \mathfrak{F} \ [\mathbf{C}_{0}^{1} \ \mathbf{C}_{1}^{1} \ \dots \ \mathbf{C}_{-2}^{1} \ \mathbf{C}_{-1}^{1}]$, the discrete Fourier transform of the first block column of~$\mathbf{C}$. By diagonalizing the block circulant matrix, the matrix-vector product~$\mathbf{T} \mathbf{p}$ can be expressed as
\begin{align}
  \mathbf{T} \mathbf{p} = 
  \tilde{\mathfrak{F}}^{-1} \lambdaMat \tilde{\mathfrak{F}} \mathbf{p}
  \text{ ,}
\end{align}
employing the modified Fourier transform~$\tilde{\mathfrak{F}}$ and its inverse~$\tilde{\mathfrak{F}}^{-1}$, i.e.,\ 
\begin{align}
  \label{eq:modifiedFourierTransform1}
  &\tilde{\mathfrak{F}} = \tilde{\mathbf{F}}_{2M_y-1} \otimes 
  \tilde{\mathbf{F}}_{2M_x-1} \otimes \mathbf{I}_{\ndof} 
  \text{ ,} \\
  \label{eq:modifiedFourierTransform2}
  &\tilde{\mathfrak{F}}^{-1} = \tilde{\mathbf{F}}_{2M_y-1}^{-1} \otimes
  \tilde{\mathbf{F}}_{2M_x-1}^{-1} \otimes \mathbf{I}_{\ndof}
  \text{ .}
\end{align}
The incomplete Fourier matrices $\tilde{\mathbf{F}}_{2M-1}$ and $\tilde{\mathbf{F}}_{2M-1}^{-1}$ contain the first $M$ columns of $\mathbf{F}_{2M-1}$ and the first $M$ rows of $\mathbf{F}_{2M-1}^{-1}$, respectively~\cite{Karimi2017}. Note that~$\lambdaMat$ is computed upfront by means of~\cref{eq:blockDiagonal} and requires the storage of~$(2M_x-1)(2M_y-1)\ndof^2$ matrix entries. This does not affect the storage costs of the PBEM since the corresponding multilevel block Toeplitz matrix can be freed from memory.

Employing this scheme, matrix-vector products with the system matrices in~\cref{eq:toeplitzBEM,eq:toeplitzBEMsymmetry} are determined in order~$\bigO(\ndof^2\Ncells\mathrm{log}(\Ncells))$ time. In the case of~\cref{eq:toeplitzBEMsymmetry}, multiplications with the first summand are performed as outlined above whereas the multiplications with the second summand, i.e.,\ with the multilevel block Hankel matrix~$\hat{\mathbf{T}}$, require a slight modification. Since Hankel matrices are column-permuted Toeplitz matrices, a permutation matrix~$\mathbf{P}$ can be applied such that~$\hat{\mathbf{T}}\mathbf{P}$ is a block Toeplitz matrix. Therefore $\hat{\mathbf{T}}\mathbf{P}\mathbf{P}^{\mathrm{T}}\mathbf{p} = \hat{\mathbf{T}}\mathbf{p}$ holds and the above scheme can also be used for multiplications with $\hat{\mathbf{T}}$ leading to the same asymptotic complexity.

\section{Fast multipole boundary element method for finite periodic arrays} \label{sec:toeplitzFMBEM}
The PBEM introduces an efficient way of assembling and storing the boundary element system matrix in the case of problems with finite periodic geometry. Regardless of the problem's size, all unique interactions between the degrees of freedom are represented by dense matrices. This inflicts unnecessary computational costs since the interaction between degrees of freedom that are well-separated can be represented in a data sparse format using fast boundary element techniques such as the fast multipole method (FMM)~\cite{Coifman1993,Darve2000}.

\subsection{Single level fast multipole method}
The fast multipole method approximates the Green's function~$G(\mathbf{x},\mathbf{y})$ by a truncated series expansion whenever the distance between a field point~$\mathbf{x}$ and a source point~$\mathbf{y}$ is sufficiently large. This decision is made based on a subdivision of the geometry into~$\nboxes$ boxes of equal size. The part of the boundary that is enclosed by a box~$\Omega_\mathbf{x}$ is assumed to be in the far field of the part of the boundary within a box~$\Omega_\mathbf{y}$ if the admissibility criterion
\begin{equation}
  \label{eq:admissibility}
  |\mathbf{x}_\mathrm{c} - \mathbf{y}_\mathrm{c}| \geq 2r
\end{equation}
holds~\cite{Darve2000}. Here, $r$ is the characteristic size of a box and $\mathbf{x}_\mathrm{c}$, $\mathbf{y}_\mathrm{c}$ are the center points of the boxes~$\Omega_\mathbf{x}$ and~$\Omega_\mathbf{y}$, respectively. Whenever~\cref{eq:admissibility} holds, the distance between a point~$\mathbf{x}\in\Omega_\mathbf{x}$ and a point~$\mathbf{y}\in\Omega_\mathbf{y}$ is sufficiently large and the Green's function can be approximated by an expansion around the center of one of the boxes. The truncated series expansion of the full-space Green's function around the center~$\mathbf{y}_\mathrm{c}$ close to~$\mathbf{y}$ reads~\cite{Liu2009}
\begin{align}
  \label{eq:greenFMM}
  G(\mathbf{x}, \mathbf{y}) \approx \dfrac{\im k}{4\pi} \sum_{n=0}^{\nt}  
    (2n + 1) \sum_{m=-n}^{n} O_n^m(\mathbf{x} - \mathbf{y}_\mathrm{c}) \bar{I}_n^m(\mathbf{y} - \mathbf{y}_\mathrm{c}) 
    \text{ ,} \qquad |\mathbf{y}-\mathbf{y}_\mathrm{c}| < |\mathbf{x}-\mathbf{y}_\mathrm{c}| \text{ ,}
\end{align}
where~$\bar{I}_n^m$ is the complex conjugate of~$I_n^m$ and~$\nt$ is the truncation number. The functions~$O_n^m$ and~$I_n^m$ are given as
\begin{align}
  \label{eq:Onm}
  O_n^m(\mathbf{x}) &= h_n^{(1)}(|\mathbf{x}|) Y_n^m \Big(\dfrac{\mathbf{x}}{|\mathbf{x}|} \Big) \text{ ,} \\
  I_n^m(\mathbf{x}) &= j_n(|\mathbf{x}|) Y_n^m \Big( \dfrac{\mathbf{x}}{|\mathbf{x}|} \Big)
  \text{ ,}
\end{align}
where $h_n^{(1)}$ denotes the $n$-th order spherical Hankel function of the first kind, $j_n$ denotes the $n$-th order spherical Bessel function of the first kind and $Y_n^m$ are the spherical harmonics~\cite{Liu2009}.

Subsequent to discretizing~\cref{eq:bm}, the integrals are split according to the admissibility criterion. Whenever the criterion holds, the truncated series expansion approximates the Green's function and allows to represent the integrals by three operators. These are the particle-to-moment (P2M), moment-to-local (M2L) and local-to-particle (L2P) operators~\cite{Liu2009}. First, the contribution of each source point~$\mathbf{y}$ within a box~$\Omega_\mathbf{y}$ is translated to the center~$\mathbf{y}_\mathrm{c}$ using the P2M operator. They are then translated to the center~$\mathbf{x}_\mathrm{c}$ of the box~$\Omega_\mathbf{x}$ by applying the M2L operator. Finally, the L2P operators translates the contribution from the center to the each field point~$\mathbf{x}$ within~$\Omega_\mathbf{x}$. Whenever the admissibility criterion does not hold, the pairs of field and source points are assumed to be in the near field and the truncated series expansion is not valid. The exact Green's function is employed and the corresponding integrals are represented by the particle-to-particle (P2P) operator. \ref{sec:appendix_fmm} presents the details of the fast multipole operators.

Assembling the discretized operators into matrices yields
\begin{align}
  \label{eq:SingleLevelFMM}
  \mathbf{H} - \mathbf{G}\mathbf{Y} \approx \Smat + \Umat \Kmat \Vmat
  \text{ .}
\end{align}
The sparse $\Ndof \times \Ndof$ matrix~$\Smat$ represents the near field interactions and the matrix product~$\Umat\Kmat\Vmat$ represents the far field interactions. $\Umat$ and $\Vmat$ are block diagonal matrices of size $\Ndof \times \nboxes(\nt + 1)^2$ and $\nboxes(\nt + 1)^2 \times \Ndof$, respectively. They store the individual L2P and P2M operators of each box. Further, $\Kmat$ is an $\nboxes(\nt + 1)^2 \times \nboxes(\nt + 1)^2$ matrix which stores the M2L operators acting between each pair of boxes that fulfills the admissibility criterion. This decomposition is known as single level fast multipole method and allows to express matrix-vector multiplications at a complexity of~$\bigO(\Ndof^{3/2})$~\cite{Coifman1993}. A further reduction can be achieved in the case of problems with periodic geometries.

\subsection{FMM for finite periodic structures}
Consider the finite periodic geometry in~\cref{fig:PeriodicArray2d} with $M_x$ and $M_y$ scatterers in the $x$ and $y$-direction, respectively. Applying the single level fast multipole method, the geometry is subdivided into boxes which correspond to the unit cell of the periodic geometry, hence~$\nboxes = \Ncells$. One of the unit cells is marked by a dashed square in~\cref{fig:PeriodicArray2d}. Similar to the PBEM, the regularity of the periodic structure in conjunction with the translation invariance of the Green's function and of its multipole expansion in~\cref{eq:multipoleExpansion} leads to the formation of special matrices. The single level FMM representation of~\cref{eq:SingleLevelFMM} reads
\begin{align}
  \label{eq:2dPeriodicFMM}
  \small
  \left(
  \underbrace{
    \begin{bmatrix}
    \Smat_{0}^{1}  		& \Smat_{-1}^{1} 	&            			&            			\\
    \Smat_{1}^{1}   	& \Smat_{0}^{1}  	& \Smat_{-1}^{1} 	&            		  \\
    & \Smat_{1}^{1}  	& \ddots 					& \ddots  				& 								\\
    &            			& \ddots 					& \ddots  				& \Smat_{-1}^{1}  \\
    &									&									& \Smat_{1}^{1} 	& \Smat_{0}^{1} 
    \end{bmatrix}
  }_{\Smat}
  +
  \underbrace{
    \mathrm{diag}(\Umat_0)
    \begin{bmatrix}
    \Kmat_{0}^{1}	& \Kmat_{-1}^{1} 	& \cdots					& \cdots & \Kmat_{1-M_y}^{1} \\
    \Kmat_{1}^{1} & \Kmat_{0}^{1}	 	& \Kmat_{-1}^{1}  & \cdots & \Kmat_{2-M_y}^{1} \\
    \vdots 				& \Kmat_{1}^{1}  	& \ddots   				& \ddots & \vdots 					 \\
    \vdots				& \vdots          & \ddots 					& \ddots & \Kmat_{-1}^{1} 	 \\
    \Kmat_{M_y-1}^{1} & \Kmat_{M_y-2}^{1} & \cdots & \Kmat_{-1}^{1} & \Kmat_{0}^{1}
    \end{bmatrix}
    \mathrm{diag}(\Vmat_0)
  }_{\Umat \Kmat \Vmat}
  \right) \mathbf{p} = \mathbf{p}^\mathrm{inc}
  \text{ .}
\end{align}
The matrix~$\Smat$ is a banded $2$-level block Toeplitz matrix which represents the near field interactions. It consists of three unique banded block Toeplitz matrices which represent the interactions between boxes within the same row $\Smat^{1}_{0}$ and in neighboring rows $\Smat^{1}_{\pm1}$. The superscript $(\cdot)^{1}$ denotes the first direction of periodicity which is the $y$-direction in the present case. In contrast, the subscript~$(\cdot)_{\pm1}$ indicates that the interaction points to the next row in positive or negative~$y$-direction, respectively. The matrix product~$\Umat\Kmat\Vmat$ represents the far field interactions. Since each box corresponds to the unit cell, the P2M and L2P operators are the same for all boxes. This leads to block diagonal matrices~$\Umat$ and~$\Vmat$ with constant blocks~$\Umat_0$ and~$\Vmat_0$, respectively. Since the translation invariance holds for the function $O^m_n(\mathbf{x})$, c.f.~\cref{eq:Onm}, it additionally holds for the M2L operator defined in~\cref{eq:m2l}. Therefore, the matrix~$\Kmat$ becomes a multilevel block Toeplitz matrix. It stores all the unique discretized M2L operators between boxes of the same row~$\Kmat^{1}_{0}$ and in different rows~$\Kmat^{1}_{j\neq0}$.

Storing the interactions of unit cells and entire rows of unit cells in this nested approach leads to a very memory efficient representation of the system matrix as long as the unit cell features a small number of degrees of freedom~$\ndof$. The presented scheme can be applied straightforwardly to geometries with periodicity in $d > 2$ directions. This yields $d$-level block Toeplitz matrices which store \mbox{$(d-1)$-level} block Toeplitz matrices itself. The inherent benefits are the same, i.e.,\ the P2M and L2P operators are equal among every unit cell, the near field matrix is a banded block Toeplitz matrix and the matrix storing the M2L operators is a multilevel block Toeplitz matrix. Similar to the PBEM, dense matrices are stored on the lowest level of the multilevel block Toeplitz matrices. In the case of~$\Smat$, they are of size~$\ndof\times\ndof$, whereas in the case of~$\Kmat$, they are of size~$(\nt+1)^2\times(\nt+1)^2$. Therefore, the total storage cost asymptotically scales~$\bigO(\ndof^2 + \nt^2 \ndof + \nt^4 M_x M_y)$ or, considering periodicity in an arbitrary number of directions, $\bigO(\ndof^2 + \nt^2 \ndof + \nt^4 \Ncells)$. This fast multipole periodic boundary element method (FMPBEM) features a more beneficial scaling than the scaling of the PBEM since $\nt \ll \ndof$ holds.

\subsection{FMM for finite periodic structures with large unit cells}
The quadratic dependence of the storage cost on~$\ndof$ limits the application to finite periodic structures with small to medium sized unit cell discretizations. A remedy is found by introducing an additional approximation to~\cref{eq:2dPeriodicFMM} that addresses the unique $\ndof\times\ndof$ matrix entries within the banded multilevel block Toeplitz matrix~$\Smat$. Applying the multilevel fast multipole method to each of the unique entries in~$\Smat$ reduces its storing cost to~$\bigO(\ndof\log(\ndof))$~\cite{Darve2000}. This leads to a nested approximation scheme where the interaction between unit cells is approximated by a single level FMM and the identified near field interactions are further reduced by a multilevel FMM.

The approximation of the unique entries of~$\Smat$ requires the subdivision of the corresponding unit cells into boxes. These are then hierarchically subdivided until a certain level of subdivision is reached or the number of degrees of freedom within a box is below a certain threshold. Boxes that are not further subdivided are called leaf boxes. Applying an admissibility criterion in the form of~\cref{eq:admissibility} allows to define the well-separated boxes on each level. The interactions are represented by P2P operators whenever the admissibility criterion does not hold. That is, the exact Green's function is employed. In contrast, interactions between well-separated boxes are represented by the truncated multipole expansion of~\cref{eq:multipoleExpansion}. Since the unit cells are hierarchically subdivided into multiple levels, the contributions are first transferred to the leaf boxes by the previously introduced P2M operator. The contributions of all leaf boxes are then transferred upwards by M2M operators, c.f.~\cref{eq:m2m}, and afterwards translated by the M2L operator. Subsequently, the contributions are transferred downwards by L2L operators, c.f.~\cref{eq:l2l}, until the leaf boxes are reached. At this stage, the contributions of the leaf boxes are transferred by the L2P operators. Introducing the truncation number~$\ntNear$ for the multipole expansion used in the multilevel FMM allows to choose its value independently of~$\nt$. Hence, the error that is added on account of the approximation of~$\Smat$ can be controlled independently from the error that is added due to the approximation of the far-field interactions between the unit cells. It is referred to~\cite{Coifman1993,Darve2000} for an in-depth analysis of the approximation error with respect to the truncation number and to~\cite{Nishimura2002,Liu2009} for details on the multilevel FMM.

The example in~\cref{fig:PeriodicArray2d} considers a structure with finite two-dimensional periodicity. Applying the FMPBEM leads to the matrix structure of~\cref{eq:2dPeriodicFMM}. Due to the quadratic shape of the unit cell in conjunction with the choice of the admissibility criterion in~\cref{eq:admissibility}, the near field matrix~$\Smat$ consists of only~$9$ unique dense matrix blocks. Each of these blocks is approximated by the multilevel FMM scheme leading to the multilevel fast multipole periodic boundary element method, referred to as FMPBEM2. The complexity of the storage cost of~$\Smat$ is reduced to~$\bigO(\ndof\log(\ndof))$ and therefore the total storage cost of the FMPBEM2 scales of order~$\bigO(\ndof\log(\ndof) + \nt^2 \ndof + \nt^4 \Ncells)$ which is quasi-linear in~$\ndof$ and linear in~$\Ncells$. The FMPBEM2 does not require to store the discretized multilevel FMM operators, however, storing them significantly accelerates the matrix-vector products with~$\Smat$. This adds an additional fourth-order dependence on the truncation number~$\ntNear$ to the asymptotically complexity which is not included here since usually~$\ntNear \ll \ndof$ holds. Whenever this does not hold, i.e.,\ the finite periodic structure contains only a small unit cell discretization, applying the FMPBEM is more favorable than applying the FMPBEM2.

\subsection{FMM for finite periodic structures in half-spaces} \label{sec:2dperiodicHalfspace}
The FMPBEM and FMPBEM2 can also be applied to acoustic half-space problems. The geometry of one half-space is modeled and the reflection of sound waves at the symmetry plane is included by employing the half-space Green's function in~\cref{eq:green3dHalf}. Whenever the symmetry plane is parallel to all directions of periodicity, the block Toeplitz matrix structure of~\cref{eq:2dPeriodicFMM} holds and a truncated series expansion of the half-space Green's function can be employed, c.f.~\cref{eq:greenFMMhalfspace} in~\ref{sec:appendix_fmm}. In all other cases, the approach outlined in~\cref{sec:CBEM_halfspace} is followed and the half-space Green's function is split into the first and the second summand. For the half-space problem shown in~\cref{fig:PeriodicArray2dSymmetry}, the system matrix is approximated by
\begin{align}
  \label{eq:2dPeriodicFMM_halfspace}
  \big((\Smat + \Umat\Kmat\Vmat) + (\hat{\Smat} + \hat{\Umat}\hat{\Kmat}\hat{\Vmat})\big) 
    \mathbf{p} 
    = 
    \mathbf{p}^{\mathrm{inc}}
  \text{ .}
\end{align}
The matrices~$\Smat$, $\Umat$, $\Kmat$ and $\Vmat$ stem from applying the fast multipole boundary element method to the first summand of the half-space Green's function. Hence, the matrices coincide with the fast multipole matrices of a full-space problem for which the symmetry plane and the mirror images are neglected. These matrices are structured matrices as shown in~\eqref{eq:2dPeriodicFMM}. The matrices~$\hat{\Smat}$, $\hat{\Umat}$, $\hat{\Kmat}$ and $\hat{\Vmat}$ stem from applying the fast multipole boundary element method to the second summand of the half-space Green's function. \ref{sec:appendix_fmm} derives the corresponding fast multipole operators. The near field matrix~$\hat{\Smat}$ is a multilevel banded block Hankel matrix whereas $\hat{\Umat}$ and $\hat{\Vmat}$ are block diagonal matrices and $\hat{\Kmat}$ is a multilevel block Hankel matrix. Note that $\mathbf{U} = \hat{\mathbf{U}}$ and $\mathbf{V} = \hat{\mathbf{V}}$ holds, since the L2P and P2M operators of both summands of the Green's function coincide. Note that the second summand includes the mirror image~$\hat{\mathbf{y}}$ instead of~$\mathbf{y}$ and thus, an additional hierarchical subdivision of the mirrored boundary is required. Consequently, the admissibility criterion for the approximation of the integrals over the second summand changes to
\begin{equation}
  \label{eq:admissibilityMirror}
  |\mathbf{x}_\mathrm{c} - \hat{\mathbf{y}}_\mathrm{c}| \geq 2r \text{ .}
\end{equation}
Herein, $\hat{\mathbf{y}}_\mathrm{c}$ is the center point of the box~$\Omega_{\hat{\mathrm{y}}}$ that encloses a part of the mirrored boundary.

\subsection{Solution scheme for the fast multipole periodic boundary element methods}
The FMPBEM and FMPBEM2 yield a system of linear equations in the form of~\cref{eq:2dPeriodicFMM} or~\cref{eq:2dPeriodicFMM_halfspace}. In either case, the matrix structure prevents the direct solution and iterative solvers have to be applied.

In the case of full-space problems, the matrix-vector operation includes multiplications with the multilevel banded block Toeplitz matrix~$\Smat$, with the block diagonal matrices~$\Umat$ and~$\Vmat$ as well as with the multilevel block Toeplitz matrix~$\Kmat$. The latter can be embedded into a multilevel block circulant matrix~$\bar{\mathbf{C}}$ which can be block diagonalized by Fourier transforms similar to the scheme outlined in~\cref{sec:PBEM_solution}. This allows to rewrite~\cref{eq:2dPeriodicFMM} into
\begin{align}
	\label{eq:FMPBEMsystem}
	(\Smat + \Umat \bar{\mathfrak{F}}^{-1} \bar{\lambdaMat} \bar{\mathfrak{F}} \Vmat)
	\mathbf{p} 
	= 
	\mathbf{p}^{\mathrm{inc}}
	\text{ ,}
\end{align}
where the block diagonal matrix~$\bar{\lambdaMat}$ stores the Fourier transform of the first block column of $\bar{\mathbf{C}}$ on its diagonal. The modified Fourier transforms~$\bar{\mathfrak{F}}$ and~$\bar{\mathfrak{F}}^{-1}$ are defined as
\begin{align}
	&\bar{\mathfrak{F}} = \tilde{\mathbf{F}}_{2M_y-1} \otimes \tilde{\mathbf{F}}_{2M_x-1} \otimes 
	\mathbf{I}_{(\nt+1)^2}		\\
	&\bar{\mathfrak{F}}^{-1} = \tilde{\mathbf{F}}_{2M_y-1}^{-1} \otimes \tilde{\mathbf{F}}_{2M_x-1}^{-1} \otimes
	\mathbf{I}_{(\nt+1)^2}
	\text{ ,}
\end{align}
with the~$(\nt+1)^2 \times (\nt+1)^2$ identity matrix~$\mathbf{I}_{(\nt+1)^2}$ and the incomplete Fourier matrices as in~\cref{eq:modifiedFourierTransform1,eq:modifiedFourierTransform2}. This reduces the asymptotic complexity of matrix-vector operations with~$\Kmat$ to~$\bigO(\nt^4\Ncells \mathrm{log}(\Ncells))$ time. Multiplications with the remaining matrices of~\cref{eq:FMPBEMsystem} are implemented as sparse matrix operations. Consequently, the computational complexity scales of order~$\bigO(\Ncells \ndof \nt^2)$ for~$\Umat$ and~$\Vmat$ as well as of order~$\bigO(\Ncells \ndof^2)$ for~$\Smat$ in the case of the FMPBEM. When employing the FMPBEM2, however, the asymptotic complexity of the multiplication with~$\Smat$ is further reduced to~$\bigO(\Ncells \ndof \log^2(\ndof))$ due to the multilevel approximation.

In the case of half-space problems for which~\cref{eq:2dPeriodicFMM_halfspace} holds, the same technique as in~\cref{sec:PBEM_solution} is applied. The matrix-vector multiplication is split into a multiplication involving the first summand and a multiplication involving the second summand. Since the first summand follows the structure of a full-space problem, i.e.,\ of~\cref{eq:2dPeriodicFMM}, this multiplication is performed as outlined above. The second summand involves block Hankel matrices and therefore a column permutation~$\mathbf{P}$ is applied to the matrix~$\hat{\Kmat}$ such that~$\hat{\Kmat}\mathbf{P}$ is a block Toeplitz matrix. This allows to address the matrix-vector product with a similar technique since $(\hat{\Smat} + \hat{\Umat}\hat{\Kmat}\hat{\Vmat})\mathbf{p} = (\hat{\Smat} + \hat{\Umat}\hat{\Kmat}\mathbf{P}\mathbf{P}^\mathrm{T}\hat{\Vmat})\mathbf{p}$ holds.

\section{Numerical examples} \label{sec:numericalExamples}
The proposed single level and multilevel fast multipole periodic boundary element methods are validated by means of two numerical examples. The first example is the scattering of a finite periodic array of spheres and the second example is a sound barrier design study. The computational efficiency of the proposed methods is compared to the conventional boundary element method (BEM), the multilevel fast multipole boundary element method (FMM,~\cite{Liu2009}) and the periodic boundary element method (PBEM). Both BEM and PBEM yield the same solution except for numerical round-off errors since the PBEM system matrix is an exact representation of the BEM system matrix. In contrast, the FMM, FMPBEM and FMPBEM2 are based on truncated multipole expansions and thus introduce additional errors. Our proposed approaches employ an optimized subdivision of the finite periodic geometry with either a single level scheme (FMPBEM) or a multilevel scheme (FMPBEM2). In contrast, the FMM uses the standard multilevel octree structure. All three methods evaluate matrix-vector products within an $l^2$-error of less than~\num{e-4} by choosing sufficiently large truncation numbers~$\nt$ as well as~$\ntNear$ in the case of the FMPBEM2.

All methods employ the generalized minimal residual method (GMRes) to solve the system of linear equations. A converged solution is found whenever a relative tolerance of~\num{e-4} is met. The calculations were performed on a desktop PC with~\SI{128}{GB} of RAM and~\num{6} physical cores running at~\SI{3.5}{GHz}. The assembly and matrix-vector operations are fully parallelized using OpenMP with $6$~threads.

\subsection{Scattering of a finite periodic array of spheres}
\begin{figure}
  \centering
  \includegraphics[width=0.5\textwidth]{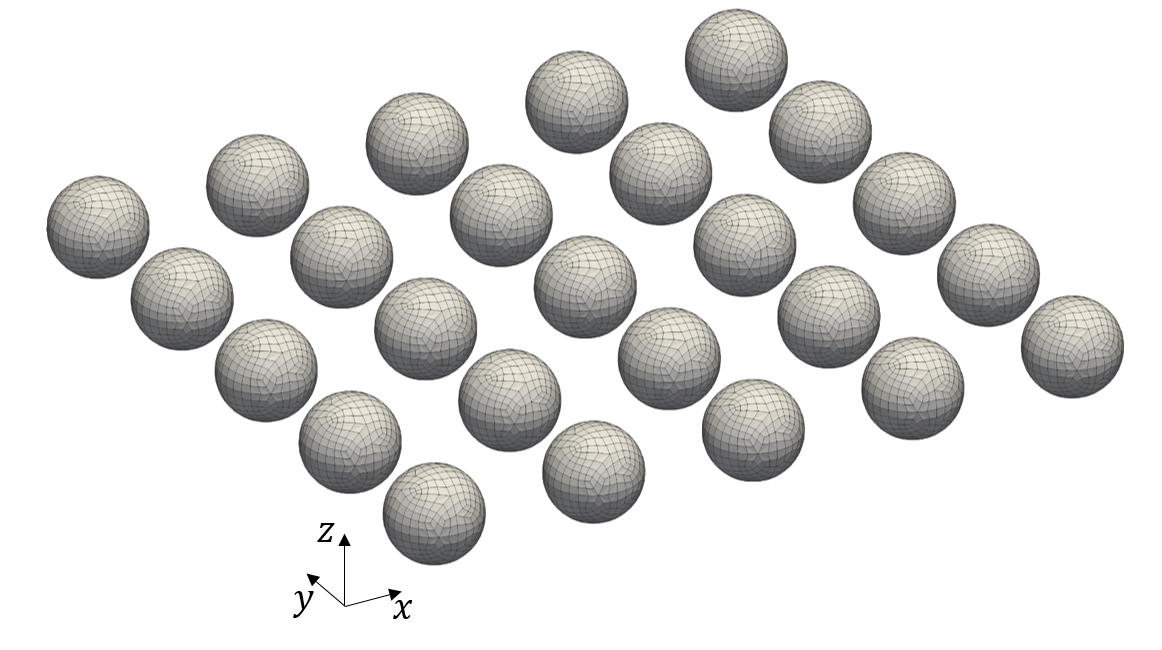}
  \caption{Finite periodic array of spherical scatters with $M_x = 5$ and $M_y = 5$. Each scatterer has a radius of~$r = \SI{100}{mm}$ and the distance between between the center points of neighboring scatterers equals $\SI{350}{mm}$ in the~$x$ and $y$-direction.}
  \label{fig:problem01_geometry}
\end{figure}

\begin{figure}
  \centering
  \includegraphics[trim=0 8 5 0, clip, width=0.75\textwidth]{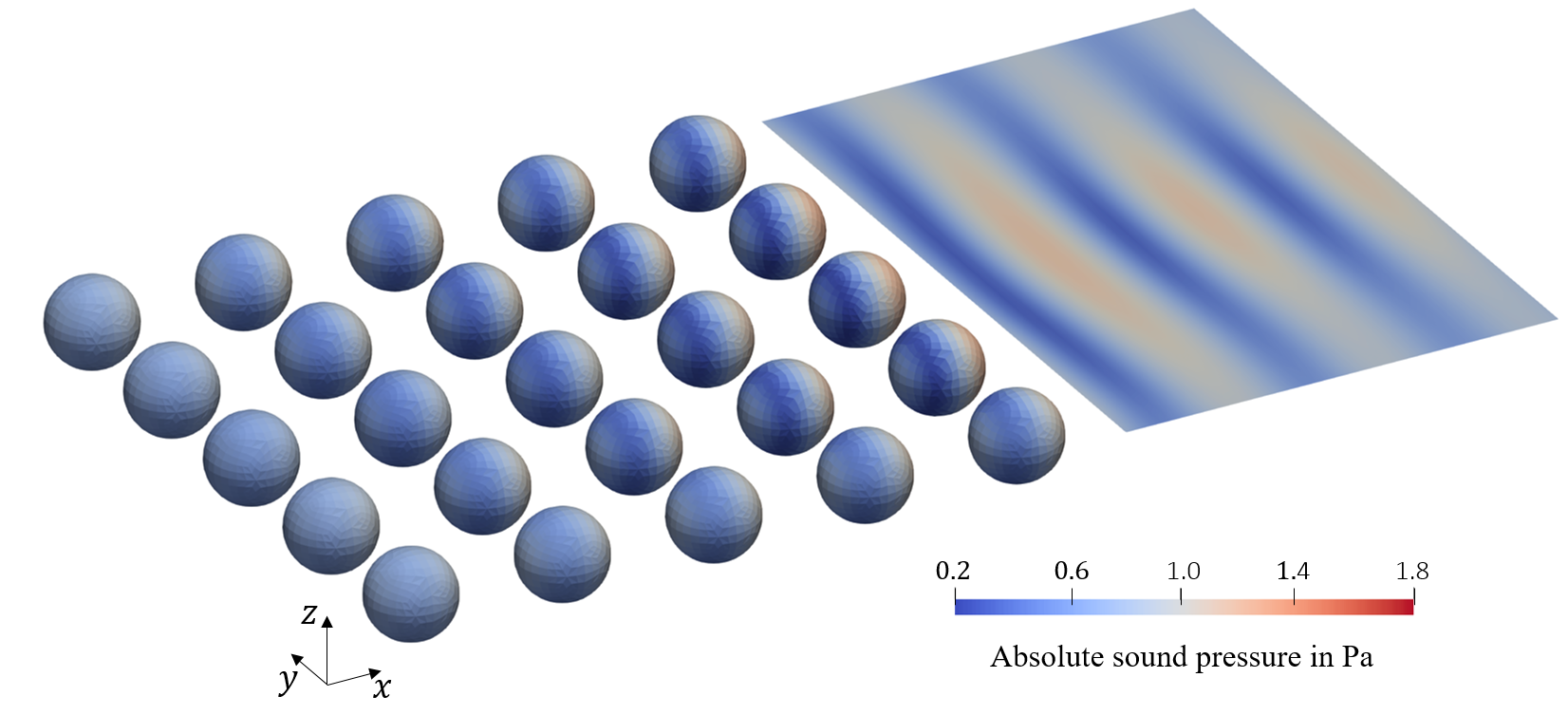}
  \caption{Absolute sound pressure on the surface of the scatterers and on a plane in the back field at~$\SI{500}{\hertz}$. Each sphere is discretized by 600 quadrilateral elements with constant discontinuous pressure approximation and excited by a plane wave traveling in the positive $x$-direction.}
  \label{fig:problem01_solution}
\end{figure}

A periodic array of $25$ acoustically rigid spheres with a radius of $r = \SI{100}{mm}$ is considered as first numerical example. The spheres are arranged in a two-dimensional pattern with $M_x = 5$ and $M_y = 5$ as shown in~\cref{fig:problem01_geometry}. The distance between the center points of neighboring scatterers in the $x$ and $y$-direction equals~$\SI{350}{mm}$. The surrounding medium is air with speed of sound of~$c = \SI{343}{\meter\per\second}$ and density of~$\rho = \SI{1,21}{\kilogram\per\meter\cubed}$. A plane wave traveling in the $x$-direction with a source strength of~$p_0=\SI{1}{\pascal}$ excites the spheres, i.e.,
\begin{align}
  \label{eq:planewave}
  p^\mathrm{inc}(x) = p_0 \ue^{\im k x}
  \text{ .}
\end{align}
Each sphere is discretized using $600$ quadrilateral boundary elements with constant discontinuous pressure approximation if not specified otherwise. This corresponds to~$40$ elements over the circumference and leads to a numerical model with~\num{15000} degrees of freedom (dofs) in total. \Cref{fig:problem01_solution} shows the absolute sound pressure on the surface of the scatterers and on a plane in the back field at~$f = \SI{500}{\hertz}$. Up to this frequency, the FMM, FMPBEM and FMPBEM2 yield solutions with a relative error of less than~\num{e-4} in the $l^2$-norm compared to the PBEM solution by prescribing truncation numbers of~$\nt = 8$ (FMM), $\nt = 4$ (FMPBEM) and $\nt = 4$, $\ntNear = 6$ (FMPBEM2).

\begin{figure}
  \begin{subfigure}[t]{0.48\textwidth}
    \centering
    \includegraphics[]{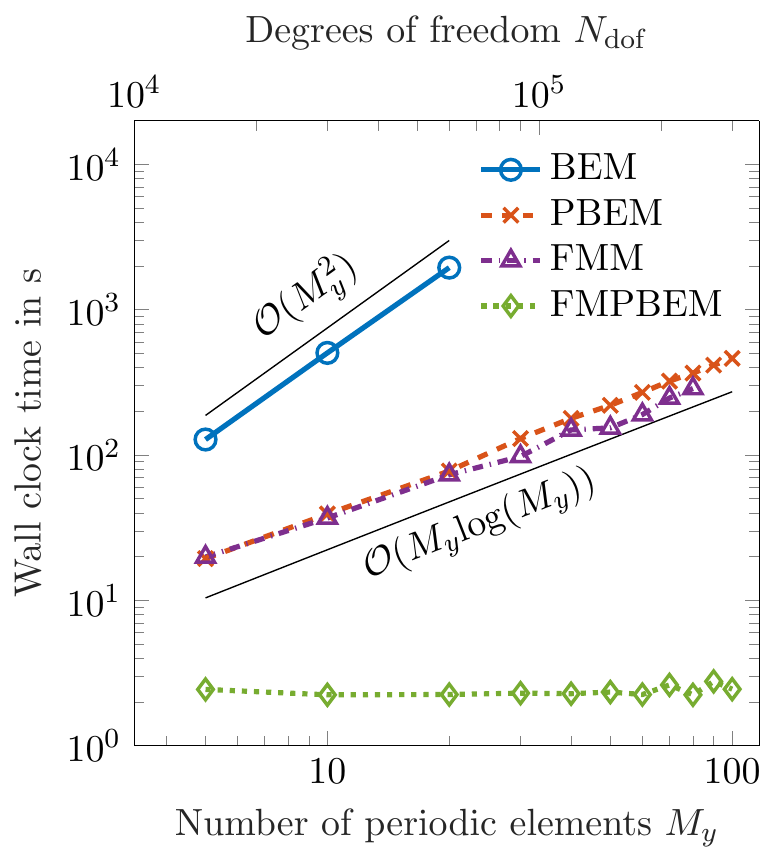}
    \vspace{-\baselineskip}
    \caption{Wall clock time of the assembly process.}
    \label{fig:problem01_timeMzAssembly}
  \end{subfigure}%
  \hfill
  \begin{subfigure}[t]{0.48\textwidth}
    \centering
    \includegraphics[]{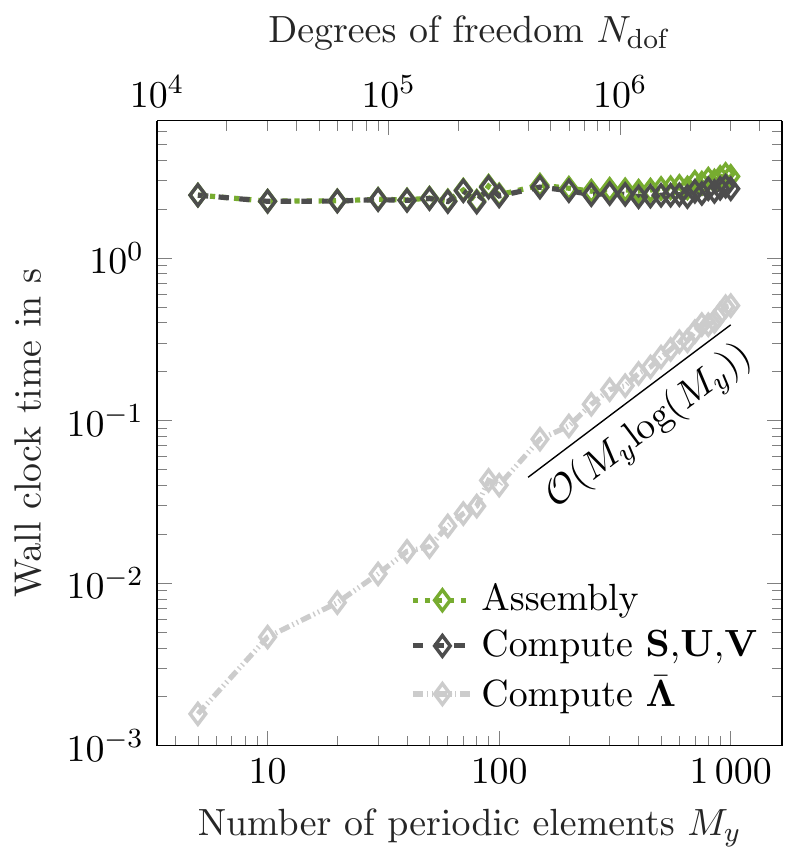}
    \caption{Wall clock time of the individual steps within the assembly of the FMPBEM. The total time adds up to the green line.}
    \label{fig:problem01_timeMzDetail}
  \end{subfigure}%
  \caption{Timings of the assembly for all five methods at~$f=\SI{500}{\hertz}$. The size of the periodic arrangement is varied by prescribing values of~$M_y$ between 5 and 100. The numerical models consist of $\Ndof= 3000M_y$ degrees of freedom.}
\end{figure}

The first study assesses the time spent for the assembly and matrix-vector computation as well as the memory usage with respect to the size of the periodic array. A variation of the number of unit cells in the $y$-direction, i.e.,\ the value of $M_y$, is performed. Values between \num{5} and \num{100} are taken into account, resulting in periodic arrangements of \num{25} to \num{500} scatterers. This corresponds to numerical models with~\num{15000} to \num{300000} dofs. \Cref{fig:problem01_timeMzAssembly} shows the time of the assembly process for the first four methods. The bottom axis represents the number of periodic elements~$M_y$ in the $y$-direction. This value relates to the total number of dofs by~$\Ndof = 3000 M_y$ which is depicted on the top axis. The BEM is applied to the first three configurations only due to its excessive memory requirements. For the initial configuration, i.e.,\ $M_y = 5$, the BEM is the slowest taking around~\SI{128.0}{\second} followed by the FMM (\SI{19.9}{\second}), the PBEM (\SI{19.4}{\second}) and the FMPBEM (\SI{2.4}{\second}). Increasing the number of unit cells in the $y$-direction, the assembly time complexity is found to be of order~$\bigO(M_y^2)$ for the BEM and of order~$\bigO(M_y \mathrm{log}(M_y))$ for the PBEM and FMM. In contrast, the assembly time of the FMPBEM seems to be constant. However, this is resolved in the more detailed plot of~\cref{fig:problem01_timeMzDetail} where values of up to~$M_y=1000$ are considered. The assembly of~$\Smat$, $\Umat$ and $\Vmat$ is constant in~$M_y$, whereas the computation of~$\bar{\lambdaMat}$ scales quasi-linearly. Therefore, the assembly of the FMPBEM matrices asymptotically scales of order~$\bigO(M_y \mathrm{log}(M_y))$. This complexity estimate equally holds for the FMPBEM2 since both methods differ only in the computation of~$\Smat$.

\begin{figure}
  \centering
  \includegraphics[]{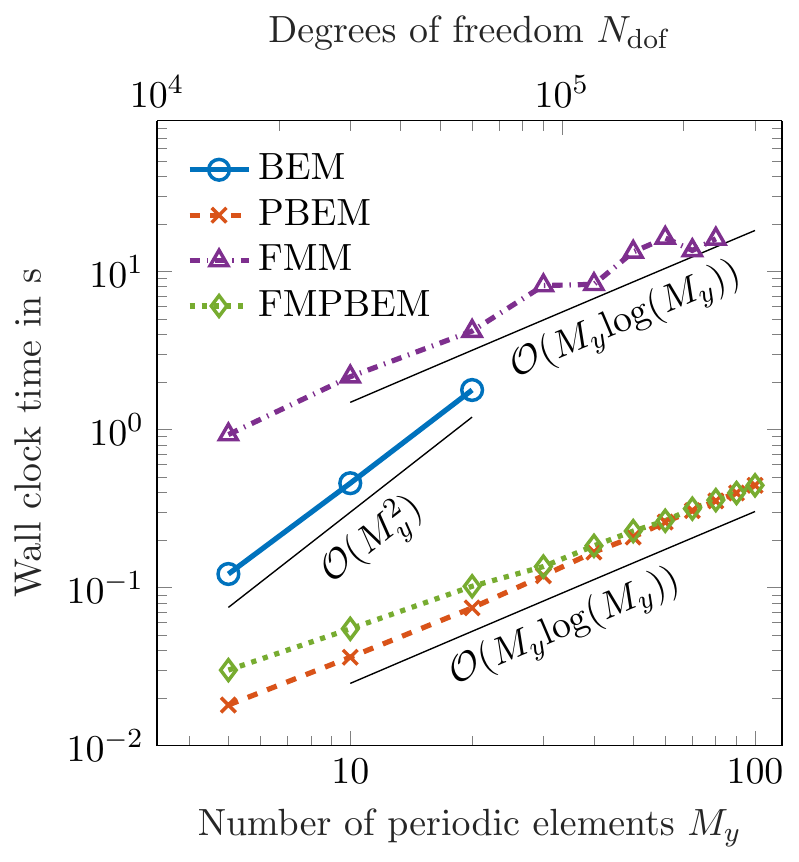}
  \caption{Wall clock time of one matrix-vector product for the first four methods at~$f=\SI{500}{\hertz}$. The size of the periodic arrangement is varied by prescribing values of~$M_y$ between 5 and 100. The numerical models consist of $\Ndof= 3000M_y$ degrees of freedom.}
  \label{fig:problem01_timeMzMatvec}
\end{figure}

\Cref{fig:problem01_timeMzMatvec} visualizes the time of one matrix-vector product for the first four methods. In the case of the initial configuration ($M_y = 5$), the FMM is the slowest taking around~\SI{0.930}{\second} followed by the BEM (\SI{0.122}{\second}), the FMPBEM (\SI{0.030}{\second}) and the PBEM (\SI{0.018}{\second}). Increasing the number of unit cells in the $y$-direction reveals the complexity of the matrix-vector products. It is of order~$\mathcal{O}(M_y^2)$ for the BEM and of order~$\mathcal{O}(M_y \mathrm{log}(M_y))$ for the PBEM and FMPBEM. The timings of the FMM follow its theoretical scaling of order~$\mathcal{O}(M_y \mathrm{log}^2(M_y))$~\cite{Darve2000} with slight fluctuations due to the changing depth of the octree subdivision with increasing~$M_y$.

\begin{figure}
	\begin{subfigure}[t]{0.48\textwidth}
	  \centering
    \includegraphics[]{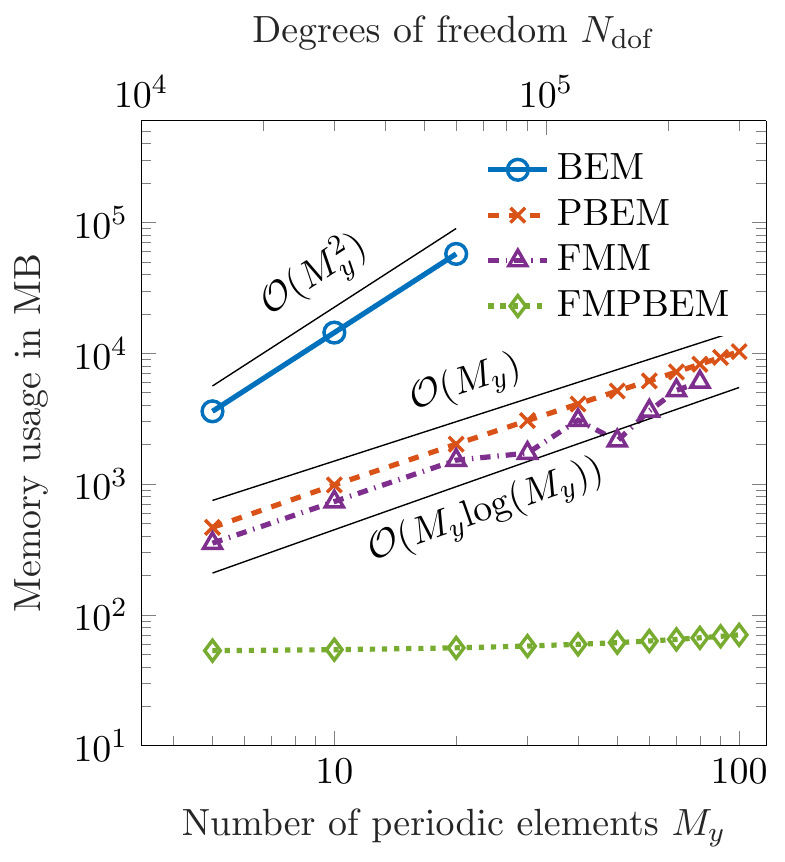}
    \vspace{-\baselineskip}
	  \caption{Allocated memory of the assembled matrices for the first four methods at $f = \SI{500}{\hertz}$. The size of the periodic array is varied by prescribing values of~$M_y$ between \num{5} and \num{100}. The numerical models consist of a total of~$\Ndof = 3000 M_y$ degrees of freedom.}
	  \label{fig:problem01_memoryMz}
	\end{subfigure}%
	\hfill
	\begin{subfigure}[t]{0.48\textwidth}
	  \centering
    \includegraphics[]{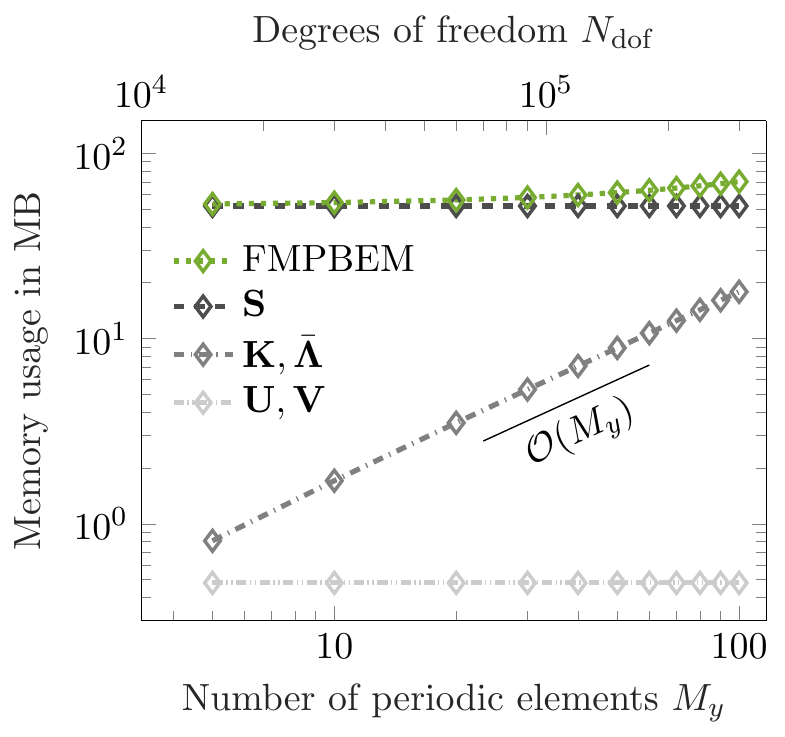}
    \vspace{-\baselineskip}
	  \caption{Allocated memory for the assembled operators of the FMPBEM. The storage costs of the matrices $\mathbf{S},\mathbf{U},\mathbf{V},\mathbf{K}$ and $\bar{\lambdaMat}$ add up to the total storage cost depicted by the green graph.}
	  \label{fig:problem01_memoryMz_detail}
  \end{subfigure}%
\caption{Comparison of the memory with respect to the number of periodic elements~$M_y$ in the $y$-direction.}
\end{figure}

\Cref{fig:problem01_memoryMz} shows the storage costs for the first four methods considering different sizes of the periodic arrangement. The BEM stores the fully populated system matrix which leads to a storage cost of~$\mathcal{O}(M_y^2)$. The memory usage of the PBEM scales linear in~$M_y$ whereas the memory usage of the FMM scales close to its theoretical order of~$\bigO(M_y \mathrm{log}(M_y))$~\cite{Darve2000}. The slight deviations can be attributed to the changing depth of the octree subdivision. \Cref{fig:problem01_memoryMz_detail} visualizes the memory usage of the FMPBEM in greater detail by illustrating the contribution of the individual matrices to the total storage costs. The allocated memory for $\mathbf{S},\mathbf{U}$ and $\mathbf{V}$ is constant in $M_y$ whereas the memory of~$\mathbf{K}$, or equivalently~$\bar{\lambdaMat}$, scales of order~$\bigO(M_y)$. Hence, the storage cost complexity of the FMPBEM asymptotically converges to~$\bigO(M_y)$.

The second study assesses the time spent for the assembly and matrix-vector computation as well as the memory usage with respect to the number of degrees of freedom~$\ndof$ within the unit cell. The layout of the periodic arrangement of scatters is held constant at~$M_x = 5$ and $M_y = 10$ whereas the boundary element discretization of each scatterer is varied. An $h$-refinement is conducted to generate the numerical models. Each sphere is initially discretized using~$24$ boundary elements with quadratic pressure approximation featuring~$192$ sound pressure dofs. The largest problem that is solved in the first comparison features~$\num{388800}$ dofs and stems from discretizing each spherical scatterer with~$972$ boundary elements with quadratic pressure approximation.

\begin{figure}
  \begin{subfigure}[t]{0.48\textwidth}
    \centering
    \includegraphics[]{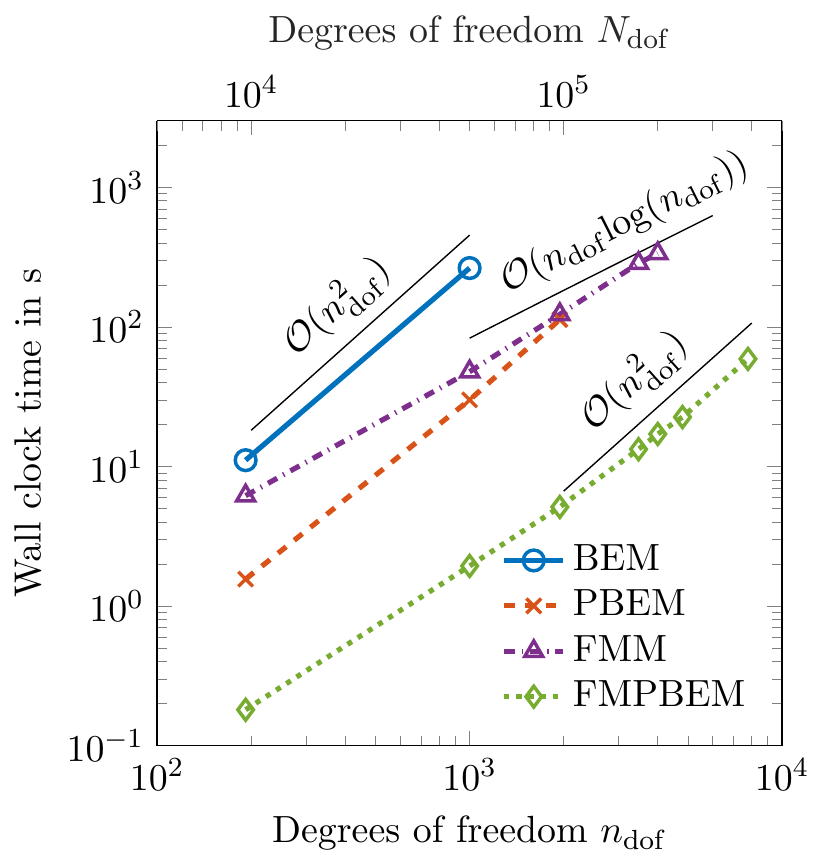}
    \vspace{-\baselineskip}
    \caption{Wall clock time of the assembly process.}
    \label{fig:problem01_time_assembly}
  \end{subfigure}%
  \hfill
  \begin{subfigure}[t]{0.48\textwidth}
    \centering
    \includegraphics[]{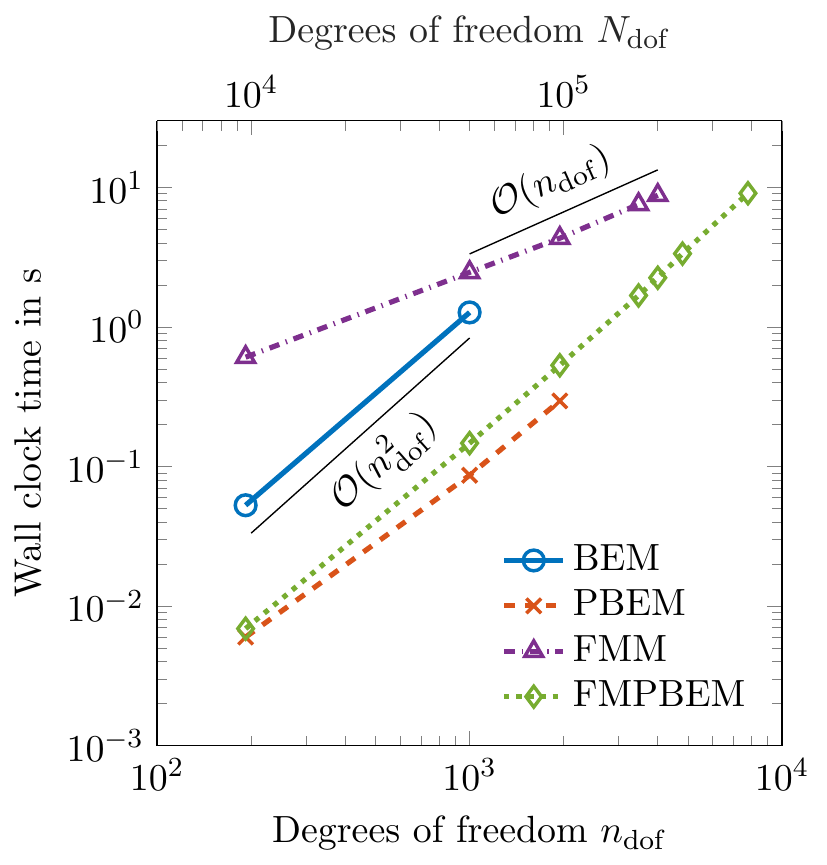}
    \vspace{-\baselineskip}
    \caption{Wall clock time of one matrix-vector product.}
    \label{fig:problem01_time_matvec}
  \end{subfigure}%
  \caption{Timings of the assembly and matrix-vector product for the first four methods at~$f=\SI{500}{\hertz}$. The number of degrees of freedom~$\ndof$ is increased by performing an~$h$-refinement. The size of the periodic arrangement is constant with $M_x = 5$ and $M_y = 10$.}
\end{figure}

\begin{figure}
  \begin{subfigure}[t]{0.48\textwidth}
    \centering
    \includegraphics[]{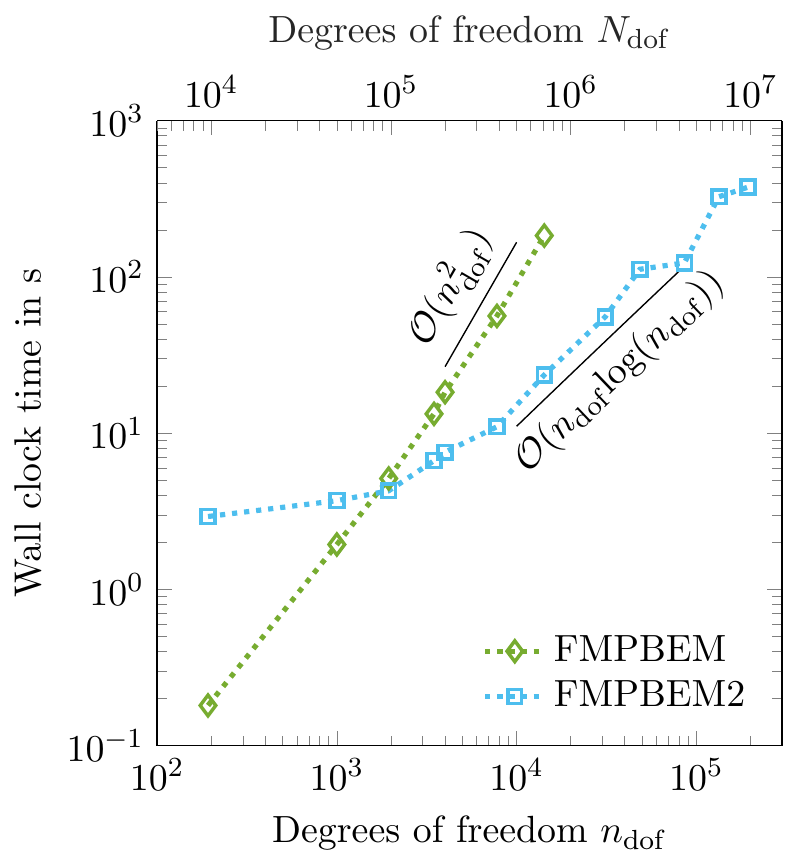}
    \vspace{-\baselineskip}
    \caption{Wall clock time of the assembly process.}
    \label{fig:problem01_time_assembly_detail}
  \end{subfigure}%
  \hfill
  \begin{subfigure}[t]{0.48\textwidth}
    \centering
    \includegraphics[]{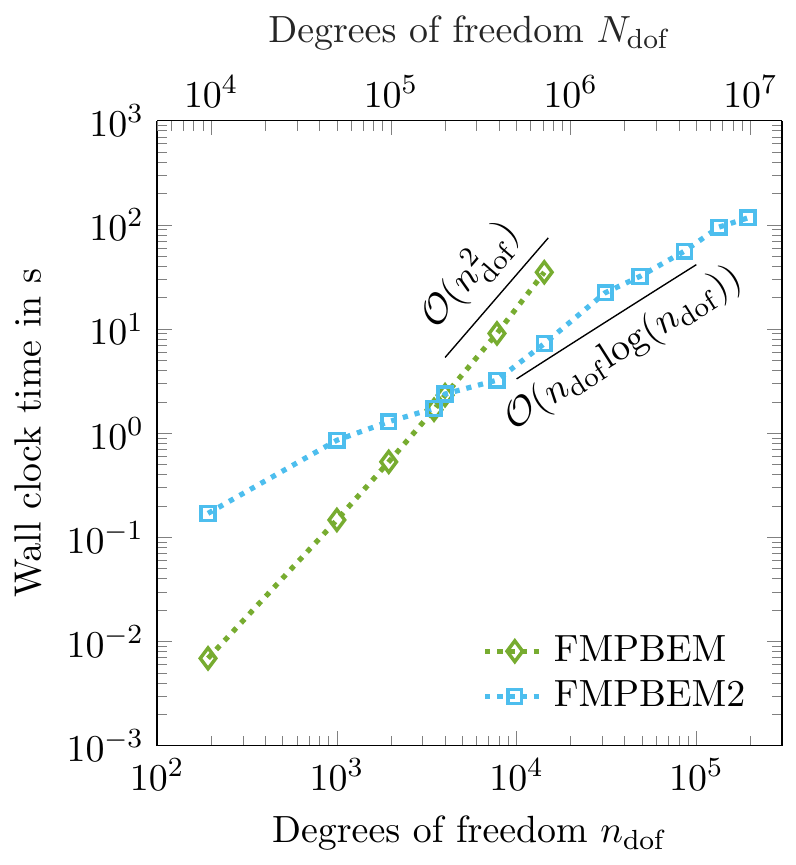}
    \vspace{-\baselineskip}
    \caption{Wall clock time of one matrix-vector product.}
    \label{fig:problem01_time_matvec_detail}
  \end{subfigure}%
  \caption{Timings of the assembly and matrix-vector product of the FMPBEM and FMPBEM2 at~$f = \SI{500}{\hertz}$ considering large-scale unit cell discretizations. The number of degrees of freedom~$\ndof$ is increased by performing an $h$-refinement while the size of the periodic arrangement is constant with~$M_x = 5$ and~$M_y = 10$.}
\end{figure}

\Cref{fig:problem01_time_assembly} visualizes the assembly time for the first four methods considering various unit cell discretizations. Not all methods are applied to every numerical model due to memory limitations. The assembly of the BEM matrices takes $\SI{11.1}{\second}$ for the smallest problem  which features a total of~$9600$ dofs. The assembly times of the FMM and PBEM are within one order of magnitude with $\SI{6.18}{\second}$ and $\SI{1.56}{\second}$, respectively. The FMPBEM assembly is the fastest, taking only $\SI{0.16}{\second}$. An increase in the number of degrees of freedom~$\ndof$ within the unit cell reveals the quadratic complexity of the BEM, PBEM and FMPBEM assembly times in~$\ndof$. In contrast, the FMM assembly time exhibits a scaling close to its theoretical value of~$\bigO(\ndof\log(\ndof))$. \Cref{fig:problem01_time_matvec} depicts the wall clock time of one matrix-vector product. The FMM takes the most time ($\SI{0.60}{\second}$) in the case of the initial configuration followed by the BEM ($\SI{0.053}{\second}$), the FMPBEM ($\SI{0.007}{\second}$) and the PBEM ($\SI{0.006}{\second}$). However, the timings of the latter three scale of order~$\bigO(\ndof^2)$ which discourages their application to models featuring very large~$\ndof$. The FMM features a scaling below the theoretical value of $\bigO(\ndof \log^2(\ndof))$ that is close to~$\bigO(\ndof)$. This might be attributed to the uneven distribution of the degrees of freedom within the~$5\times{}10\times{}1$ pattern of the spherical scatterers. Employing a multilevel approximation of the near field matrix within the FMPBEM reduces the computational complexity of the assembly and matrix-vector computation. \Cref{fig:problem01_time_assembly_detail} and~\cref{fig:problem01_time_matvec_detail} visualize the performances of the FMPBEM in comparison to its extension, the FMPBEM2. For small unit cell discretizations, the FMPBEM is faster in both assembly and matrix-vector operations due to the additional overhead of the approximation of the near field matrix~$\Smat$ in the FMPBEM2. However, in the case of medium to large-scale unit cell discretizations, the FMPBEM2 achieves a significant reduction in computational time due to its more favorable scaling of order~$\bigO(\ndof\log(\ndof))$. The largest problem features~$\ndof=\num{194144}$ degrees of freedom within each unit cell leading to a total of~$\Ndof=\num{9707200}$ dofs. Its assembly takes $\SI{376}{\second}$ and one matrix-vector product is computed in~$\SI{117}{\second}$. The same problem cannot be solved with the FMM due to memory limitations.

\begin{figure}
  \begin{subfigure}[t]{0.48\textwidth}
    \centering
    \includegraphics[]{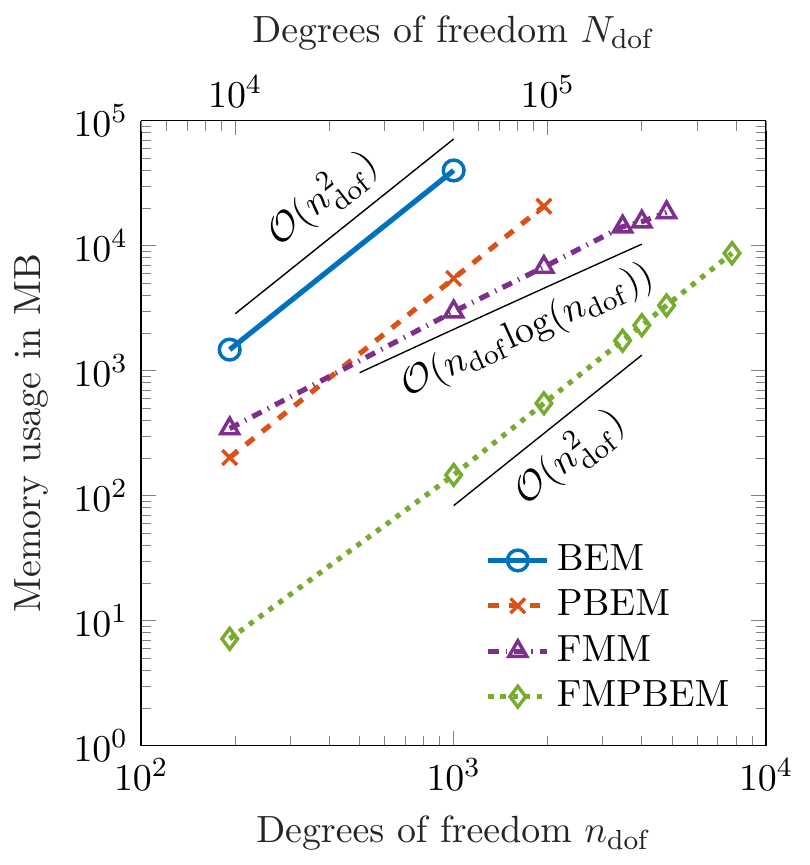}
    \vspace{-\baselineskip}
    \caption{Storage cost of the assembled matrices for the first four methods.}
    \label{fig:problem01_memory}
  \end{subfigure}%
  \hfill
  \begin{subfigure}[t]{0.48\textwidth}
    \centering
    \includegraphics[]{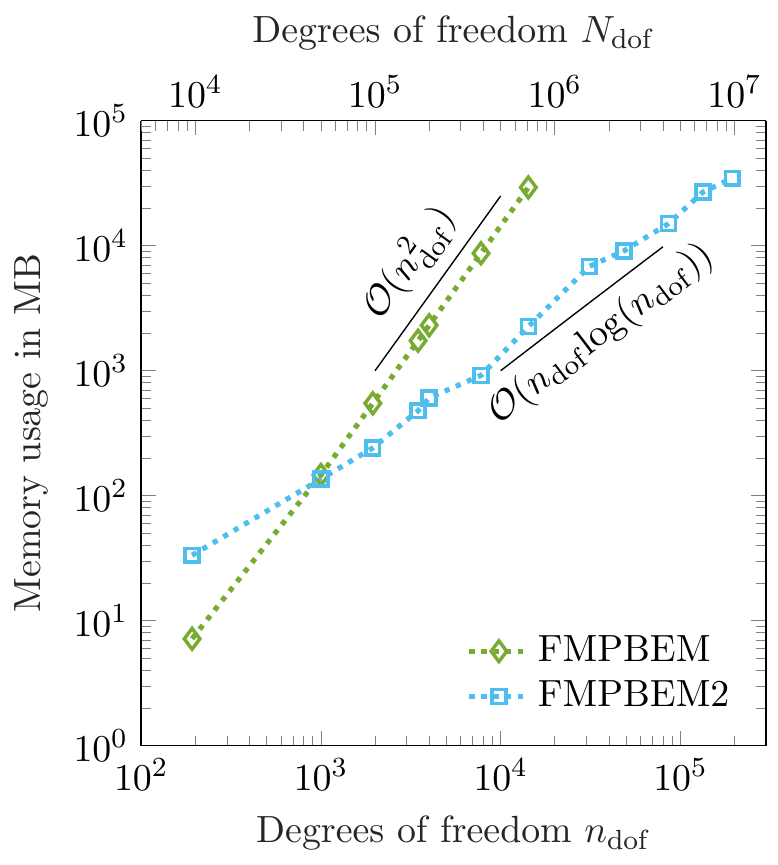}
    \vspace{-\baselineskip}
    \caption{Storage cost of the assembled matrices for the FMPBEM and FMPBEM2 considering large-scale unit cell discretizations.}
    \label{fig:problem01_memory_detail}
  \end{subfigure}%
  \caption{Comparison of the storage cost with respect to the number of degrees of freedom~$\ndof$ at~$f = \SI{500}{\hertz}$. The size of the periodic arrangement is constant with $M_x = 5$ and $M_y = 10$.}
\end{figure}

\Cref{fig:problem01_memory} shows the memory usage of the first four methods considering varying unit cell discretizations. The BEM, PBEM and FMPBEM exhibit a quadratic scaling of the storage cost with respect to the number of degrees of freedom~$\ndof$ of the unit cell. In contrast, the FMM scales of order~$\bigO(\ndof \mathrm{log}(\ndof))$. Although the FMPBEM features the lowest memory usage within the considered range of discretizations, the quadratic scaling prevents its application to finite periodic problems with large values of~$\ndof$. This is resolved by the FMPBEM which introduces an additional approximation of the near field matrix~$\Smat$. \Cref{fig:problem01_memory_detail} visualizes the memory usage for large-scale unit cell discretizations with up to~$\ndof=\num{194114}$ degrees of freedom and indicates a complexity of~$\bigO(\ndof\log(\ndof))$ for the storage cost within the FMPBEM2.

\begin{figure}
  \centering
  \includegraphics[]{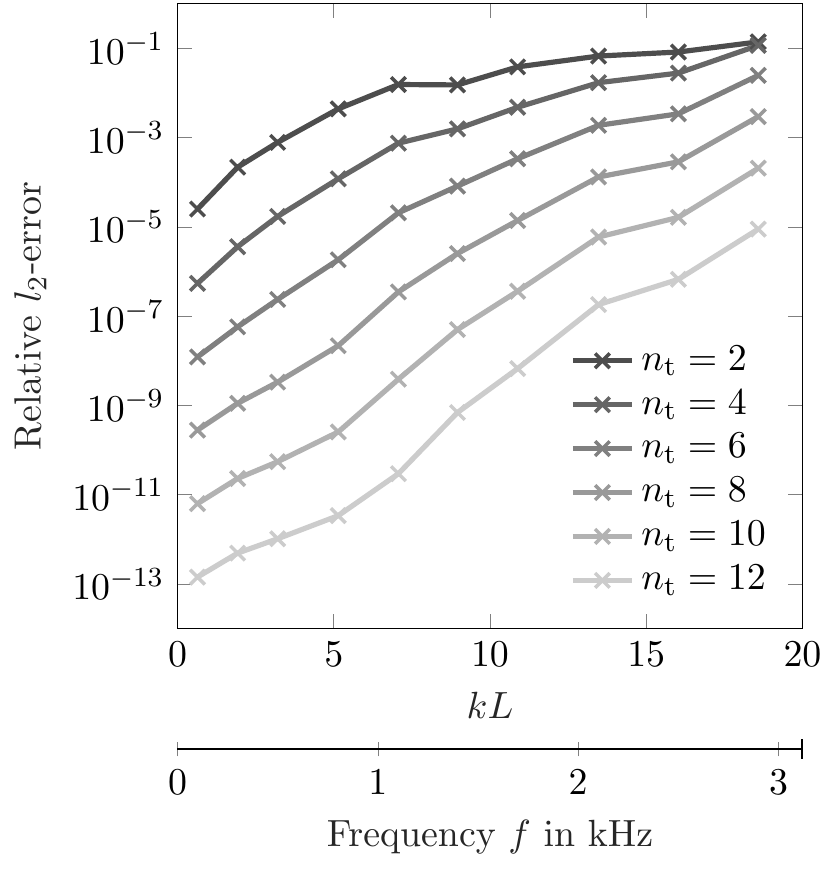}
  \caption{Relative error of the FMPBEM solution to the BEM solution in the~$l_2$-norm as a function of the dimensionless wavenumber $kL$, with $L$ denoting the characteristic unit cell length. A periodic arrangement with $M_x = 5$ and $M_y = 100$ is analyzed. Each sphere is discretized with~$600$ boundary elements with constant pressure approximation. Different values of the truncation number~$n_{\mathrm{t}}$ are considered.}
  \label{fig:problem01_error}
\end{figure}

The third study analyzes the frequency dependent accuracy of the FMPBEM solution. The study includes a periodic arrangement with~$M_x = 1$ and~$M_y  = 30$ and frequencies of up to~$\SI{3000}{\hertz}$. Each sphere is discretized using $600$ boundary elements with quadratic pressure approximation. This corresponds to about~$7.3$ elements per wavelength at~$\SI{3}{\kilo\hertz}$. The GMRes tolerance is set to~$\num{e-14}$ and truncation numbers of~$\nt = 2{,}\,4{,}\,\ldots{,}\,12$ are considered. \Cref{fig:problem01_error} presents the relative error in the $l_2$-norm of the FMPBEM solution to a reference solution determined by the PBEM. The $x$-axis shows both, the frequency $f$ and the dimensionless wavenumber $kL$ with the characteristic length~$L = \SI{350}{\milli\meter}$ of the unit cells. A relative error of less than~$\num{e-4}$ is achieved for all considered truncation numbers at~\SI{100}{\hertz} or $kL = \num{0.64}$. With an increase in frequency, the accuracy of the FMPBEM deteriorates. A relative error of less than~$\num{e-4}$ at~$\SI{1000}{\hertz}$ requires six or more terms of the fast multipole expansion.

\subsection{Sound barrier}
\begin{figure}
  \centering
  \begin{subfigure}[t]{0.40\textwidth}
    \centering
    \includegraphics[width=0.95\textwidth]{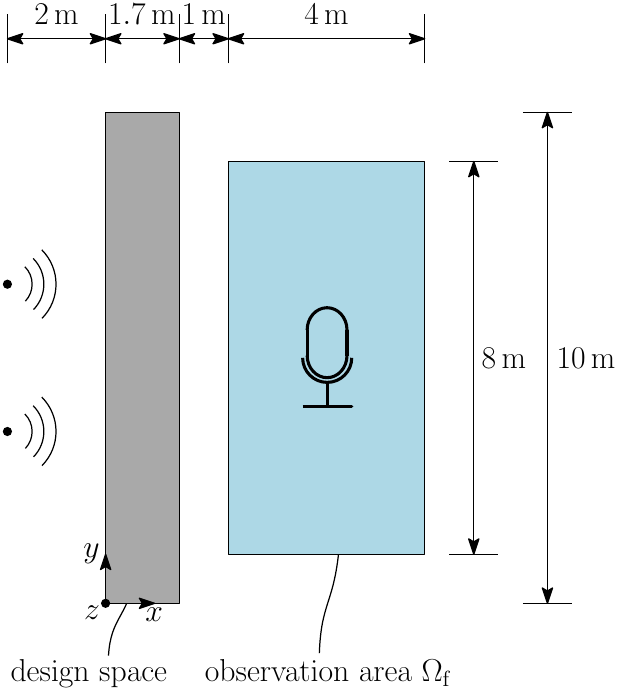}
    \caption{}
    \label{fig:problem02_setup}
  \end{subfigure}%
  \hfill
  \begin{subfigure}[t]{0.60\textwidth}
    \centering
    \adjincludegraphics[Clip={.12\width} {.1\height} {0.2\width} {.05\height},width=0.8\textwidth]{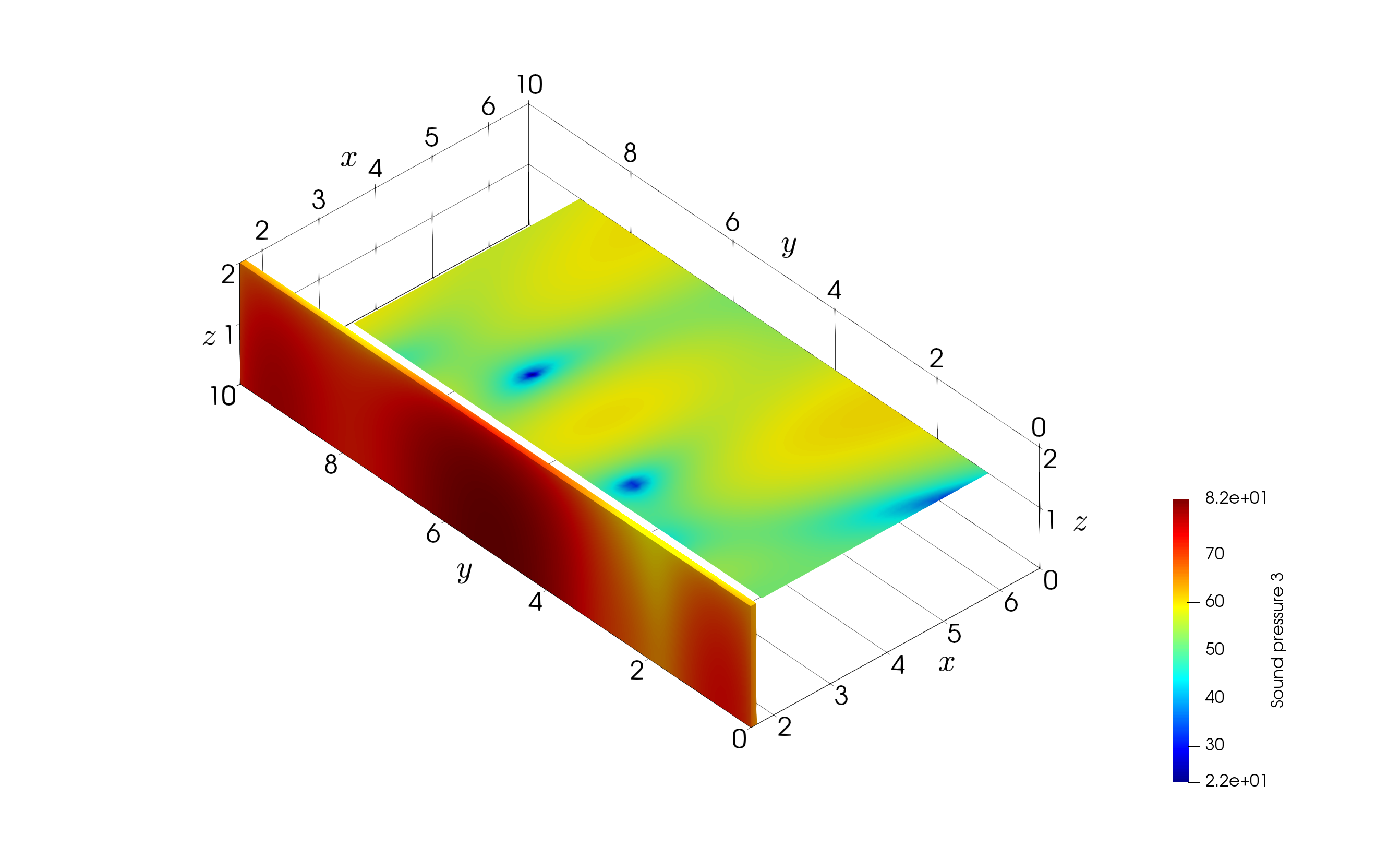}
    \includegraphics[]{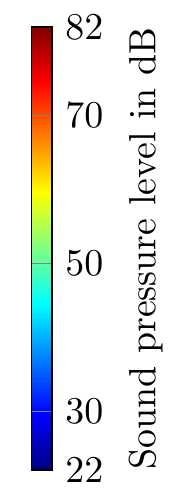}
    \vspace{-\baselineskip}
    \caption{}
    \label{fig:problem02_wall}
  \end{subfigure}
  \caption{Top view of the sound barrier setup (a) and sound pressure level ($\SI{0}{\decibel} = \SI{2e-5}{\pascal}$) on the wall sound barrier and in the observation area at~$\SI{100}{\hertz}$ (b).}
\end{figure}

A sound barrier design study is considered in the second numerical example. The general setup is shown in~\cref{fig:problem02_setup}. Two monopole sources emit an incident sound pressure wave on the left-hand side, $\SI{1}{\meter}$ above the sound-hard ground. They differ in the source strength which is $p_1 = \SI{2}{\pascal}$ for the source at~$y_1=\SI{6.5}{\meter}$ and $p_2 = \SI{1}{\pascal}$ for the source at~$y_2=\SI{3.5}{\meter}$. A sound barrier is located within the design space and ideally reduces the sound pressure within the observation area~$\Omega_\mathrm{f}$. Assessing the insertion loss (IL) within the observation area quantifies the performance of different sound barrier designs. Following the work of~\cite{Cavalieri2019}, the IL reads
\begin{equation}
  \label{eq:insertionloss}
  \mathrm{IL}(f) = 
  20 \log_{10} \left(\frac{\sum_{i=1}^{M} |\mathbf{p}^{\mathrm{inc}}(\mathbf{x}_i,f)|} {\sum_{i=1}^{M} |\mathbf{p}(\mathbf{x}_i,f)|}\right) \mathrm{ ,} \quad
  \mathbf{x}_i \in \Omega_{\mathrm{f}}
  \text{ ,}
\end{equation}
with the number of observation points~$M$ within~$\Omega_f$. In this study, $M=3200$ uniformly distributed points are taken into account. A frequency range of~$\SI{100}{\hertz}$ to $\SI{500}{\hertz}$ is analyzed with $81$ uniformly distributed frequency samples.

\begin{figure}
  \centering
  \begin{subfigure}[t]{0.3\textwidth}
    \centering
    \includegraphics[width=0.85\textwidth]{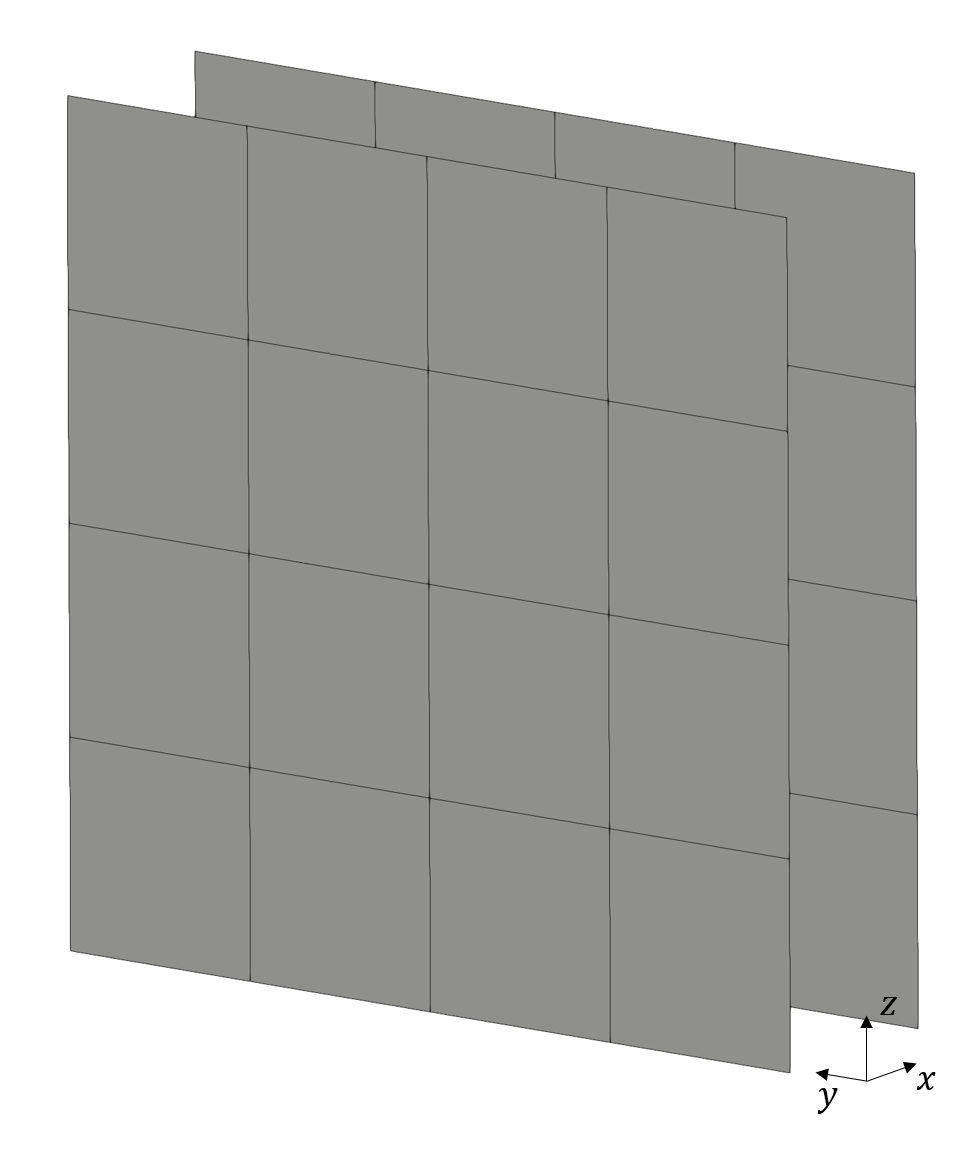}
    \caption{$\SI{0.2}{\meter} \times \SI{0.2}{\meter} \times \SI{0.1}{\meter}$\\ $M_x = 1$, $M_y = 50$, $M_z = 10$}
    \label{fig:problem02_unitcell_wall}
  \end{subfigure}%
  \hfill
  \begin{subfigure}[t]{0.3\textwidth}
    \centering
    \includegraphics[width=0.55\textwidth]{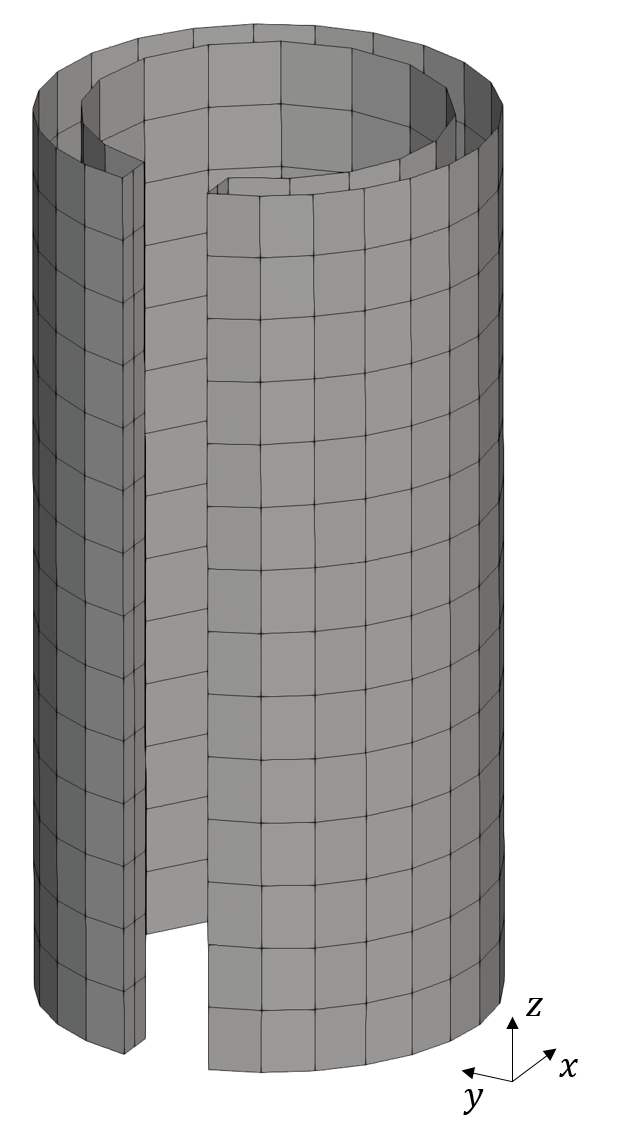}
    \caption{$\SI{0.2}{\meter} \times \SI{0.2}{\meter} \times \SI{0.4}{\meter}$\\ $M_x = 3$, $M_y = 25$, $M_z = 5$}
    \label{fig:problem02_unitcell_cylinder}
  \end{subfigure}%
  \hfill
  \begin{subfigure}[t]{0.3\textwidth}
    \centering
    \includegraphics[width=0.55\textwidth]{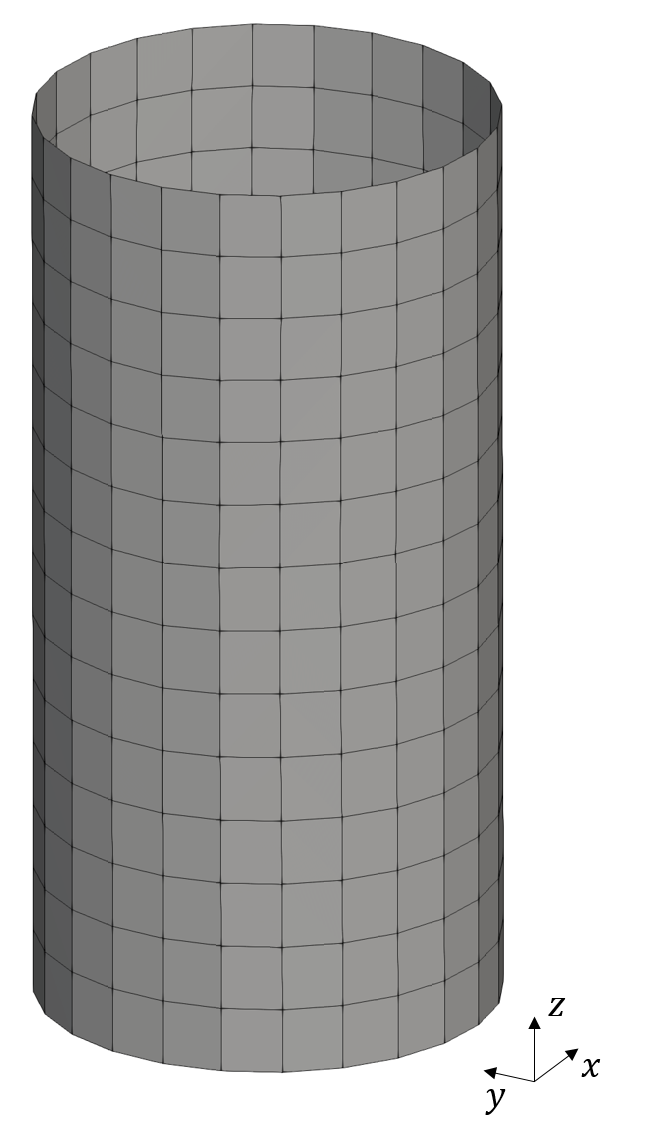}
    \caption{$\SI{0.2}{\meter} \times \SI{0.2}{\meter} \times \SI{0.4}{\meter}$\\ $M_x = 3$, $M_y = 25$, $M_z = 5$}
    \label{fig:problem02_unitcell_cshape}
  \end{subfigure}%
  \caption{Unit cells of the wall sound barrier (a), cylinder sound barrier (b) and c-shaped sound barrier (c). The outer dimensions of the unit cells and the number of cells in each direction of periodicity are given below the corresponding figures.}
\end{figure}

The first sound barrier design is a sound-hard wall with constant rectangular cross section, height of~$\SI{2}{\meter}$, length of~$\SI{10}{\meter}$ and a width of~$\SI{0.1}{\meter}$. It is located at the rightmost part of the design space. The full-scale wall model is shown in~\cref{fig:problem02_wall} and consists of $4140$ boundary elements with quadratic pressure approximation. This equals~$\num{6.9}$ elements per wavelength at~$\SI{500}{\hertz}$ and a total of~$33120$ pressure dofs. The unit cell of the corresponding periodic wall model is shown in~\cref{fig:problem02_unitcell_wall}. It consists of two boundary element layers, one at the front and one at the back. The layers have a size of~$\SI{0.2}{\meter}$ by~$\SI{0.2}{\meter}$ and consist of four boundary elements with quadratic pressure approximation each. This also leads to~$\num{6.9}$ elements per wavelength at~$\SI{500}{\hertz}$. The unit cell is extended by~$M_y = 50$ cells in the $y$-direction, i.e.,\ the length of the wall, and $M_z = 10$ cells in the $z$-direction, i.e.,\ the height of the wall. Hence, a total of $500$ periodic cells are considered. Note that the periodic model lacks the top and side surfaces of the wall and therefore only consists of~$32000$ dofs.

\begin{figure}
  \centering
  \includegraphics[]{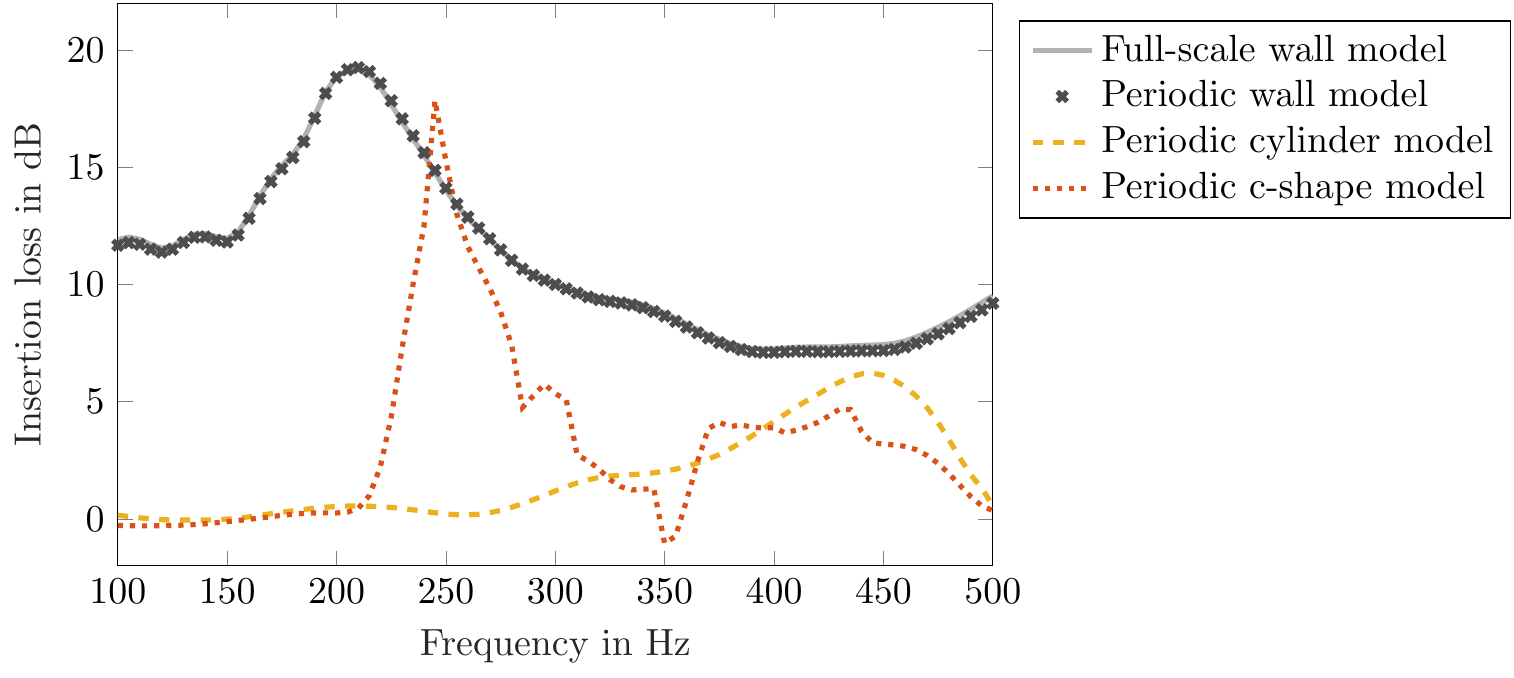}
  \caption{Insertion loss (IL) in \si{\decibel} of the three sound barrier designs. The IL is assessed in the observation area~$\Omega_{\mathbf{f}}$ according to~\cref{eq:insertionloss}.}
  \label{fig:problem02_insertionloss_monopole}
\end{figure}

\Cref{fig:problem02_insertionloss_monopole} shows the insertion loss of the sound-hard wall over the frequency range of~$\SI{100}{\hertz}$ to~$\SI{500}{\hertz}$ for both models, the full-scale wall model and the periodic wall model. The solution of the former is generated by the BEM whereas the FMPBEM is employed for solving the latter. Although the periodic model lacks the top and side surfaces of the wall, the insertion loss values are in good agreement with the results of the full-scale model and underestimate the IL only slightly. The insertion loss stays well above~$\SI{10}{\decibel}$ up to~$\SI{300}{\hertz}$ with a peak value of~$\SI{19.2}{\decibel}$ at~$\SI{210}{\hertz}$. The minimum insertion loss value of~$\SI{7.6}{\decibel}$ is found at~$\SI{390}{\hertz}$. \Cref{fig:problem02_wall} shows the sound pressure levels on the sound barrier surface and in the observation area for~$f=\SI{100}{\hertz}$.

The second sound barrier design considers a periodic array of cylinders where each cylinder has a radius of~$\SI{0.1}{\meter}$ and a height of~$\SI{2}{\meter}$. A total of $75$ cylinders are arranged in a rectangular pattern with~$M_x = 3$ rows and~$M_y = 25$ cylinders along the length of the barrier. With a distance of~$d = \SI{0.4}{\meter}$ between neighboring cylinders in the~$x$ and $y$-direction, the footprint of the cylinder sound barrier wall is~$\SI{9.8}{\meter}$ by~$\SI{1}{\meter}$. The unit cell has a height of~$\SI{0.4}{\meter}$, c.f.~\cref{fig:problem02_unitcell_cylinder}, and is repeated~$M_z = 5$ times in the~$z$-direction, resulting in a total height of~$\SI{2}{\meter}$. Each unit cell is discretized with $12$ boundary elements with linear sound pressure approximation along the circumference and~$35$ elements along its height. This corresponds to at least~$12$ elements per wavelength at~$\SI{500}{\hertz}$ and results in~$126000$ dofs in total.

Although the cylinder sound barrier is not continuous along its length, a significant insertion loss is still expected due to standing waves forming between the cylinders. According to Bragg's law, standing waves occur at 
\begin{equation}
  f_\mathrm{Bragg} = \frac{n c}{2 d \, \mathrm{sin}(\theta)}
  \text{ ,}
\end{equation}
with a positive integer~$n$ and the incident angle~$\theta$. In the present case, the first Bragg frequency equals $\SI{428,75}{\hertz}$ assuming perpendicular incidence, i.e.,\ $\theta=\pi/2$. Note that~$0.125\pi < \theta \leq \pi/2$ holds for the given monopole excitation and sound barrier setup. The insertion loss within the frequency range of interest is depicted in~\cref{fig:problem02_insertionloss_monopole}. The IL stays around~$\SI{0}{\decibel}$ up to~$\SI{280}{\hertz}$ from which on it monotonically increases to its maximum of~$\SI{6.5}{\decibel}$ at~$\SI{440}{\hertz}$.

\begin{figure}
  \begin{subfigure}[t]{0.42\textwidth}
    \centering
    \adjincludegraphics[Clip={.1\width} {.08\height} {0.2\width} {.04\height},width=0.95\textwidth]{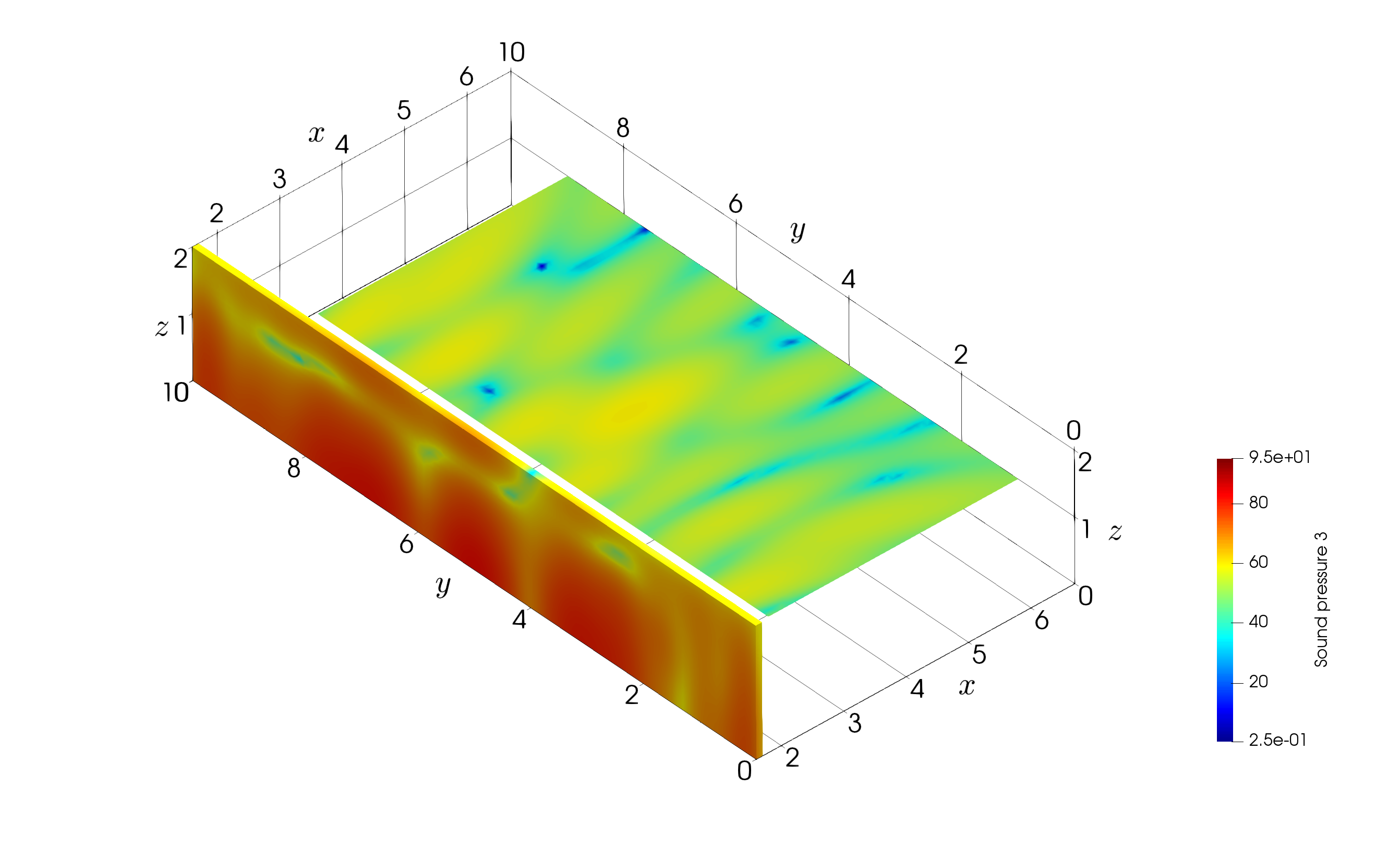}
    \caption{Wall sound barrier}
  \end{subfigure}%
  \begin{subfigure}[t]{0.42\textwidth}
    \centering
    \adjincludegraphics[Clip={.07\width} {.05\height} {0.2\width} {.04\height},width=0.95\textwidth]{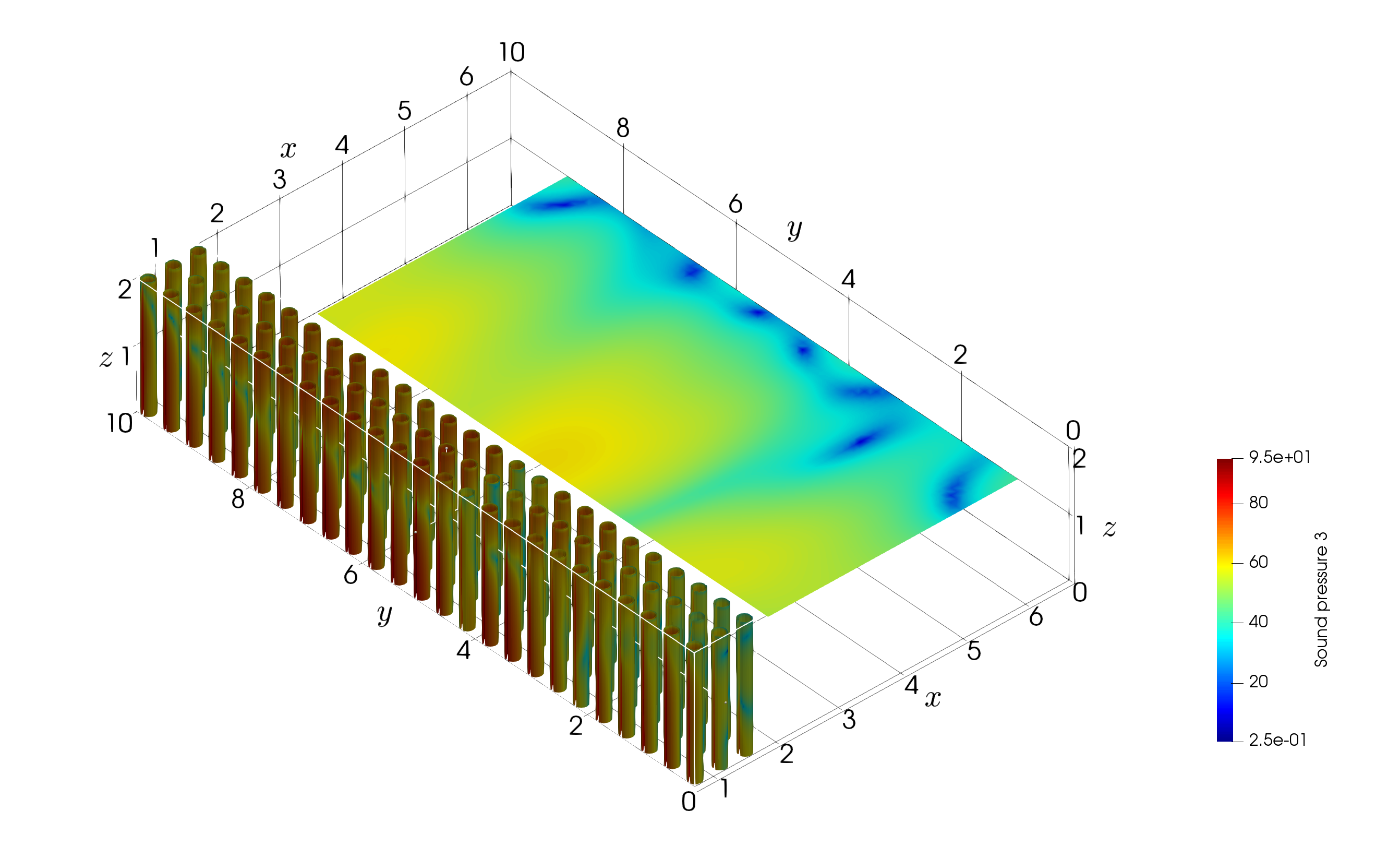}
    \caption{C-shaped sound barrier}
  \end{subfigure}%
  \begin{subfigure}[t]{0.16\textwidth}
    \includegraphics[]{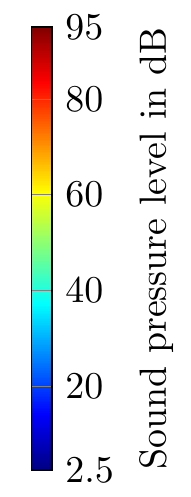}
  \end{subfigure}%
  \hfill
  \caption{Sound pressure level ($\SI{0}{\decibel} = \SI{2e-5}{\pascal}$) on two sound barrier designs and in the observation area at~$\SI{245}{\hertz}$.}
  \label{fig:problem02_solution}
\end{figure}

The third sound barrier design follows the idea of the second design but introduces an additional resonance by considering a c-shaped cross section. Each c-shape structure has an outer radius of~$\SI{0.1}{\meter}$, inner radius of~$\SI{0.08}{\meter}$, a slit width of~$\SI{0.04}{\meter}$ and a height of~$\SI{2}{\meter}$. The unit cell of the periodic model has a height of~$\SI{0.4}{\meter}$ and is depicted in~\cref{fig:problem02_unitcell_cshape}. Its discretization features around $11.9$ boundary elements with linear pressure approximation per wavelength at~$\SI{500}{\hertz}$. Besides the standing waves around the Bragg frequency an additional resonance occurs since each c-shape structure acts as a Helmholtz resonator. The resonance frequency is estimated at~$f_\mathrm{HR}=\SI{258.9}{\hertz}$ based on a two-dimensional finite element simulation of the c-shaped cross section. The parameters of the periodic array remain unchanged with $M_x = 3$, $M_y = 25$, $M_z=5$ and~$d = \SI{0.4}{\meter}$ which leads to the same estimate of the Bragg frequency~$f_\mathrm{Bragg}$ and a total of~$231000$ dofs. The IL is shown in~\cref{fig:problem02_insertionloss_monopole} and proceeds very similar to the IL of the cylinder sound barrier design up to around~$\SI{210}{\hertz}$. From here on, the IL of the c-shaped design rises steeply to a maximum of~$\SI{17.9}{\decibel}$ at~$\SI{245}{\hertz}$. This is followed by an equally steep decline up to~$\SI{320}{\hertz}$. A further significant IL is reported between~$\SI{370}{\hertz}$ and~$\SI{440}{\hertz}$ with values above~$\SI{3.8}{\decibel}$. Introducing a c-shaped cross section significantly improves the insertion loss within a frequency region around~$f_\mathrm{HR}$ while still utilizing the Bragg effects. \Cref{fig:problem02_solution} shows the sound pressure values of the wall sound barrier and the c-shaped sound barrier in a side-by-side comparison at a frequency of~$f = \SI{245}{\hertz}$.

\section{Conclusion and future work}
Two fast multipole periodic boundary element methods have been proposed for the solution of time-harmonic acoustic problems of finite periodic structures. Both methods subdivide the geometry into boxes that correspond to the unit cells of the periodic geometry. A boundary element discretization is applied to each unit cell, i.e.,\ each box. Interactions between well-separated boxes are approximated by a multipole expansion of the Green's function. On the other hand, either an exact representation of the Green's function is used for the interactions of neighboring boxes (FMPBEM) or an additional multilevel fast multipole method is employed (FMPBEM2). In both methods, the fast multipole operators acting between the unit cells become block Toeplitz matrices due to the periodicity of the geometry. These are matrices with constant blocks along each diagonal and therefore feature only a small memory footprint. In addition, matrix-vector products can be expressed by circular convolutions which significantly reduces their computational complexity. Certain configurations of half-space problems disturb the block Toeplitz structure. However, this contribution introduces a remedy by splitting the half-space Green's function into it's summands and discretizing each summand separately. The first numerical example, the scattering of a periodic array of sound-hard spheres, indicated a quasi-linear scaling of the fast multipole periodic boundary element solution with respect to the number of unit cells. The solution time has been found to be an order of magnitude below comparable approaches such as the periodic boundary element method and the multilevel fast multipole method. The study indicated that the proposed methods perform especially well in cases where large periodic structures are considered. The application of the FMPBEM is favorable in the case of unit cells that feature small numbers of degrees of freedom, whereas the FMPBEM2 is more beneficial in the case of large-scale unit cell discretizations. The accuracy of our approach is set by the truncation number of the multipole expansions which also has been investigated in the first example. In the second example, the FMPBEM has been applied to a sound barrier design study in half-space. Both, wall sound barriers as well as sonic crystal sound barriers were analyzed. Although the periodic model of the wall barrier does not account for the top and side surfaces, the results were in good agreement with the results of a full-scale analysis. It has been shown that a comparison of sound barrier designs in three dimensional space is indeed feasible with the proposed methods.

The concept of the fast multipole periodic boundary element methods is based on the translation invariance of the Green's function. Applying the same concept to problems with a rotationally arranged unit cells would yield a method that, for instance, would be able to analyze the aeroacoustic behavior of ducted fans more efficiently. Furthermore, the presented approach can be applied likewise using the hierarchical boundary element method. With respect to the computational efficiency, further improvements are planned by representing the multilevel block Toeplitz matrices in the tensor train format, c.f.~\cite{Polimeridis2014} for details. Future work will involve taking visco-thermal losses into account and introducing a structural-acoustic coupling scheme. This will allow to analyze sound barriers that include acoustic energy harvesting as presented in~\cite{Wang2018}. In addition, the proposed methods can be extended to study local defects within the periodic structure by introducing cost-efficient low-rank updates of $\Smat$, $\Umat$ and $\Vmat$.

\section*{Acknowledgments}
\addcontentsline{toc}{section}{Acknowledgements}
This work is financially supported by the China Scholarship Council (CSC) (File No. 201706340085), the National Natural Science Foundation of China (NSFC) under Grant No. 11772322.

\appendix
\renewcommand{\thesection}{Appendix \Alph{section}}
\section{Fast multipole expansions}
\label{sec:appendix_fmm}
A single integral of~\cref{eq:bm} over an arbitrary boundary part~$\Gamma_\mathrm{c} \subset \Gamma$ is picked as an example. All source points~$\mathbf{y}$ on the boundary~$\Gamma_\mathrm{c}$ lie within the box~$\Omega_{\mathbf{y}}$. This box is in the far field of the box~$\Omega_{\mathbf{x}}$ which encloses the field points $\mathbf{x}$. Based on the truncated series expansion of the full-space Green's function~\cref{eq:greenFMM}, the approximation of the integral reads~\cite{Liu2009}
\begin{equation}
	\begin{split}
  \label{eq:multipoleExpansion}
  \int_{\Gamma_\mathrm{c}} 
    \dfrac{\partial G(\mathbf{x}, \mathbf{y})}{\partial \mathbf{n}(\mathbf{y})}
    p(\mathbf{y}) \diff \Gamma(\mathbf{y})
  \approx
  \ &\dfrac{\im k}{4\pi} \sum_{n=0}^{\nt} (2n + 1) \sum_{m=-n}^{n} 
    O_n^m(\mathbf{x} - \mathbf{y}_\mathrm{c}) M_{n}^{m}(\mathbf{y}, \mathbf{y}_\mathrm{c})
  \text{ ,} \\
    &|\mathbf{y}-\mathbf{y}_\mathrm{c}| < |\mathbf{x}-\mathbf{y}_\mathrm{c}| 
  \text{ ,}
  \end{split}
\end{equation}
with the multipole moments~$M_{n}^{m}$ at the expansion point~$\mathbf{y}_\mathrm{c}$ given as
\begin{align}
  \label{eq:p2m}
  M_{n}^{m}(\mathbf{y}_\mathrm{c})
  = 
  \int_{\Gamma_\mathrm{c}} \dfrac{\partial \bar{I}_{n}^{m}(\mathbf{y} -
    \mathbf{y}_\mathrm{c})}{\partial \mathbf{n}(\mathbf{y})}
  p(\mathbf{y}) \diff \Gamma(\mathbf{y})
  \text{ .}
\end{align}
\Cref{eq:multipoleExpansion} is called multipole expansion and~\cref{eq:p2m} is the particle-to-multipole (P2M) translation. The discretization of the latter equation leads the P2M operator.

Expanding the Green's function around a point~$\mathbf{x}_\mathrm{c}$ close to~$\mathbf{x}$ instead yields the local expansion~\cite{Liu2009}
\begin{equation}
	\begin{split}
	\label{eq:l2p}
  \int_{\Gamma_\mathrm{c}} 
    \dfrac{\partial G(\mathbf{x}, \mathbf{y})}{\partial \mathbf{n}(\mathbf{y})}
    p(\mathbf{y}) \diff \Gamma(\mathbf{y})
  \approx
  \ &\dfrac{\im k}{4\pi} \sum_{n=0}^{\nt} (2n + 1) \sum_{m=-n}^{n} 
    \bar{I}_n^m(\mathbf{x} - \mathbf{x}_\mathrm{c})
    L_{n}^{m}(\mathbf{x}_\mathrm{c}) 
  \text{ ,} \\
    &|\mathbf{x} - \mathbf{x}_\mathrm{c}| < |\mathbf{y} - \mathbf{x}_\mathrm{c}|
  \text{ ,}
  \end{split}
\end{equation} 
with the local coefficients~$L_{n}^{m}$ at the expansion point~$\mathbf{x}_\mathrm{c}$ given as
\begin{align}
  \label{eq:m2l}
  L_{n}^{m}(\mathbf{x}_\mathrm{c}) 
  =
  \sum_{n'=0}^{\nt} (2n' + 1) \sum_{m'=-n'}^{n'} (-1)^{m+m'} \sum_{l \in N} 
    W_{n',n,m',m,l} O_{l}^{m+m'}(\mathbf{x}_\mathrm{c} - \mathbf{y}_\mathrm{c}) 	
    M_{n'}^{m'}(\mathbf{y}_\mathrm{c})
  \text{ .}
\end{align}
\Cref{eq:l2p} is the local-to-particle (L2P) translation and~\cref{eq:m2l} is the multipole-to-local (M2L) translation. The set~$N$ is defined by~\cite{Liu2009}
\begin{align}
  N := \Big\{l \ | \ l \in \mathbb{Z},\, n + n' - l: \mathrm{even},\, 
  \mathrm{max} \left\{|m + m'|, |n - n'| \right\} < l < n + n'\Big\}
  \text{ ,}
\end{align}
and $W_{n',n,m',m,l}$ is given as
\begin{align}
  \label{eq:wigner}
  W_{n',n,m',m,l} = (2l + 1) \im^{n'-n+l} 
  \begin{pmatrix}
    n & n' & l \\
    0 & 0  & 0 
  \end{pmatrix}
  \begin{pmatrix}
    n & n' & l \\
    m & m' & -m-m' 
  \end{pmatrix}
  \text{ ,}
\end{align}
where~$(:\,:\,:)$ denotes the Wigner 3j-symbol~\cite{Shore1968}. 

In the case of half-space problems, a truncated series expansion is employed for the second summand of the half-space Green's function~\cref{eq:green3dHalf}. The expansion reads
\begin{equation}
  \label{eq:greenFMMhalfspace}
  G(\mathbf{x}, \hat{\mathbf{y}}) \approx 
    R_\mathrm{p} \dfrac{\im k}{4\pi} \sum_{n=0}^{\nt} (2n + 1) \sum_{m=-n}^{n} 
    O_n^m(\mathbf{x} - \hat{\mathbf{y}}_\mathrm{c}) 
    \bar{I}_n^m(\hat{\mathbf{y}} - \hat{\mathbf{y}}_\mathrm{c}) 
  \text{ ,} \qquad |\hat{\mathbf{y}}-\hat{\mathbf{y}}_\mathrm{c}| 
    < |\mathbf{x}-\hat{\mathbf{y}}_\mathrm{c}| \text{ .}
\end{equation}
The mirrored source points~$\hat{\mathbf{y}}$ lie in the box~$\Omega_{\hat{\mathbf{y}}}$ with center point~$\hat{\mathbf{y}}_\mathrm{c}$. \Cref{eq:greenFMMhalfspace} is valid whenever the admissibility criterion holds, i.e., when $\Omega_{\hat{\mathbf{y}}}$ is in the far field of $\Omega_{\mathbf{x}}$. The corresponding fast multipole operators can be derived similarly to the aforementioned operators of the full-space problem by substituting~$\mathbf{y}$ and~$\mathbf{y}_\mathrm{c}$ with its mirrored variants.

In the case of the multilevel fast multipole method, two additional operators are introduced. The moment-to-moment (M2M) translation shift the multipole moments from an expansion point~$\mathbf{y}_\mathrm{c}$ to an expansion point~$\mathbf{y}_{\mathrm{c}\prime}$ following
\begin{equation}
  \label{eq:m2m}
    \tilde{M}_{n}^{m}(\mathbf{y}_{\mathrm{c}\prime}) 
    = \sum_{n'=0}^{\nt} (2n' + 1) \sum_{m'=-n'}^{n'} 
    \sum_{l \in N} (-1)^{m'} W_{n,n',m,m',l}
    I_{l}^{-m-m'}(\mathbf{y}_\mathrm{c} - \mathbf{y}_{\mathrm{c}\prime}) 	
    M_{n'}^{-m'}(\mathbf{y}_\mathrm{c})
  \text{ .}
\end{equation}
Similarly, the local-to-local (L2L) translation shifts the local coefficients from an expansion point~$\mathbf{x}_{\mathrm{c}\prime}$ to an expansion point~$\mathbf{x}_{\mathrm{c}\prime}$ by
\begin{equation}
  \label{eq:l2l}
    L_{n}^{m}(\mathbf{x}_{\mathrm{c}\prime}) 
    = (-1)^m \sum_{n'=0}^{\nt} (2n' + 1) \sum_{m'=-n'}^{n'} \sum_{l \in N} 
    W_{n',n,m',m,l} 
    I_{l}^{m-m'}(\mathbf{x}_{\mathrm{c}\prime} - \mathbf{x}_\mathrm{c})
    L_{n'}^{m'}(\mathbf{x}_\mathrm{c})
  \text{ .}
\end{equation}

\bibliographystyle{elsarticle-num-names}
\bibliography{tfmmPaper}

\end{document}